\documentclass{tac}
\usepackage[utf8]{inputenc}
\usepackage{amsfonts}
\usepackage{hyperref}
\usepackage{amsmath}
\usepackage{bbm}
\usepackage{braket}
\usepackage{tikz-cd}
\usepackage{enumitem}
\usepackage{mathtools}
\usepackage{hyphenat}
\usepackage{tabularx}
\usepackage{bm}
\usepackage{mathpartir}

\makeatletter
\newcommand{\pto}{}
\newcommand{\pgets}{}

\DeclareRobustCommand{\pto}{\mathrel{\mathpalette\p@to@gets\to}}
\DeclareRobustCommand{\pgets}{\mathrel{\mathpalette\p@to@gets\gets}}
\newcommand{\p@to@gets}[2]{%
  \ooalign{\hidewidth$\m@th#1\mapstochar\mkern5mu$\hidewidth\cr$\m@th#1\to$\cr}%
}
\makeatother

\title[Commutativity and Kleisli laws of codensity monads]{Commutativity and Kleisli laws for codensity monads of probability measures}

\author{Zev Shirazi}
\copyrightyear{2024}
\begin{document}
\maketitle

\begin{abstract}
Several monads of probability measures have been shown to have presentations as codensity monads over small categories of stochastic maps. This paper studies how three key properties of these probability monads, relevant to categorical approaches to probability, can arise from their codensity presentations. We first derive the existence of a Kleisli law into the Giry monad, which provides a formal connection to measurable probability. In particular, from their codensity presentations, we prove a novel universal property of several probability monads as terminal liftings of the Giry monad. This generalises a result by Van Breugel on the Kantorovich monad, and proves the existence of such Kleisli laws. We additionally provide sufficient conditions for a codensity monad to be lax monoidal and affine, which provides a connection to the theory of Markov categories. In particular, we introduce the condition for a codensity monad to be exactly pointwise monoidal, which is then lax monoidal, and prove a characterisation of this condition in terms of Day convolution. We show that the Radon monad is exactly pointwise monoidal, and use our characterisation to give a description of the tensor product of free algebras of the Radon monad in terms of Day convolution. Finally, we show that the Giry monad is only exactly pointwise monoidal when restricted to standard Borel spaces,  due to the existence of probability bimeasures that do not extend to measures. 
\end{abstract}

\tableofcontents
\newpage
\section{Introduction}
Monads have risen to become an important and versatile tool for probability theory, due to the variety of probabilistic examples and results that can be described from a monadic perspective ~\cite{giry_categorical_1982,fritz_finettis_2021}. Monads for probability have been used in semantics of probabilistic programming languages~\cite{jones_probabilistic_1989,heunen_convenient_2017}, and functional analysis~\cite{semadeni_monads_1973,tadeusz_swirszcz_monadic_1974}. In the past ten years, several new \emph{probability monads} have been studied~\cite{fritz_probability_2019,jacobs_probability_2018} and have seen applications in synthetic approaches to probability theory via \emph{Markov categories}~\cite{fritz_synthetic_2020,cho_disintegration_2019}. A formal categorical perspective on these monads has been developed, which shows that spaces of probability measures are \emph{categorical limits}, and that monads for probability can be derived from this limit structure, called \emph{codensity monads}~ \cite{avery_codensity_2016,van_belle_probability_2022}. This paper builds on the work of \cite{van_belle_probability_2022} by discussing when \emph{properties of probability monads} can additionally be derived from this limit structure. We focus on three properties of probability monads which are relevant to synthetic approaches to probability \cite{fritz_synthetic_2020,jacobs_probability_2018}.
\begin{enumerate}
       \item  We provide basic conditions for codensity monads to be \emph{affine} in \hyperref[prop3.11]{Proposition 3.11}. This provides part of a condition for the Kleisli category of a codensity monad to be a Markov category.
    \item In \hyperref[section4]{Section 4}, we show the existence of a Kleisli law~\cite{mulry_lifting_1993}, or monad morphism, into the Giry monad \cite{giry_categorical_1982}. This connects the codensity monad presentation of probability~\cite{van_belle_probability_2022} to the historical measure-theoretic presentation of probability~\cite{halmos_lectures_1974}. This is shown by developing a connection between codensity presentations of probability monads and a novel universal property as \textit{liftings} of the Giry monad. This has been previously studied for the Kantorovich monad~\cite{van_breugel_metric_nodate}, but we generalise the extent of this result in \hyperref[thm4.9]{Theorem 4.9} and provide similar characterisations for other probability monads. 
    \item It is well-known that monads with a \emph{lax monoidal structure}~\cite{kock_monads_1970,kock_strong_1972} are equivalent to \emph{commutative monads}~\cite{kock_closed_1971}. In a programming language setting, commutative monads model order-invariant effect usage. Probability provides an example of a commutative effect, and in the measure-theoretic setting, this amounts to Fubini's theorem~\cite{rudin_functional_1991}. Furthermore, combined with affineness, this provides a condition for the Klesli category of a probability monad to form a Markov category. In \hyperref[section5]{Section 5}, we show when codensity monads have a lax monoidal structure, which allows for analysis of bivariate probability distributions. In particular, we introduce the notion of an \textit{exactly pointwise monoidal codensity monad}, which is an analogue of codensity monads in the setting of monoidal categories, and is related to the theory of algebraic Kan extensions~\cite{fritz_criterion_nodate,weber_algebraic_2016}. We characterise this notion in \hyperref[thm5.23]{Theorem 5.23}, which is the main result of this section.
\end{enumerate}
\subsection{Background on probability monads}  

If $A$ and $B$ are finite sets, a random function $f \colon A \pto  B$, which given an input $a \in A$ takes a value $b \in B$ with some probability $f(b|a)\in [0,1]$, can be reinterpreted as an ordinary function $\tilde{f} \colon A \to \mathcal{D}_cB$ where $ \mathcal{D}_cB$ is the set of discrete probability measures on $B$. Here, $\tilde{f}(a)$ is taken to be the probability distribution on $B$, defined by $(f(b|a))_{b \in B}$. More generally, if $A$ and $B$ are measurable spaces, a Markov kernel $f \colon A \pto B$ which, given an input $a \in A$ takes a value in a measurable subset $E\subseteq B$ with a probability $f(E|a)\in [0,1]$, can be reinterpreted as a measurable map $\tilde{f} \colon A \to \mathcal{G}B$ where $\mathcal{G}B$ is the set of probability measures on $B$ with a suitable $\sigma${-}algebra.
Composition of random functions, given by
\[g \bullet f(c|a)= \sum_{b\in B}\: g(c|b)\cdot f(b|a)\]
and of Markov kernels, given by the Chapman-Kolmogorov formula~\cite{lawvere_category_1962}
\[{g \bullet f}(E | a)=\int_B g(E|{-}) \; \text{d}f({-}|a)\]
suggest these random functions can be understood as the Kleisli morphisms of a monad. There are numerous examples of similar \emph{probability monads}~\cite{giry_categorical_1982,jacobs_probability_2018}, and we will describe five running examples for this paper in \hyperref[section2]{Section 2}. These monads, although defined on different categories, tend to share a similar description in terms of integration operators, shown in \hyperref[table1]{Table 1}.

\begin{table} [htbp]
    \centering
    \begin{tabularx}{0.85\textwidth}{|p{1cm}|p{5cm}<{\centering}|X<{\centering}|} \hline 
         $PX$&  Space of probability measures on $X$& Initial for maps   $\int_X h \;  \text{d}{-} \colon PX \to [0,1]$ for $h\in \mathcal{C}(X,[0,1])$\\ \hline 
         $Pf$&  Pushforward of measures along $f \colon X \to Y$& Unique map such that $\int_Y h \;  \text{d}Pf(p) = \int_X hf \; \text{d}p$\\ \hline 
         $\eta_X$&  Dirac measure of $x\in X$ & Unique map such that $\int_X h 
 \; \text{d}\eta_X(x) =h(x)$\\ \hline 
         $\mu_X$&  Average value of a measure $\mathbb{P}$ on probability measures in $PX$ & Unique map such that $\int_{PX} h  \; \text{d}\mu_X(\mathbb{P}) = \int (\int h_X \; \text{d}{-}) \;  \text{d}\mathbb{P}$ \\ \hline \end{tabularx}
\caption{Description of a probability monad $P$ on $\mathcal{C}$}
\label{table1}
\end{table}
This unified description can be explained by a presentation of probability monads as codensity monads~\cite{van_belle_probability_2022}. If $K \colon \mathcal{D} \to \mathcal{C}$ is any functor, then if the right Kan extension $T=\text{Ran}_K K$ exists, it has a monad structure, which is called the codensity monad of $K$. In particular, if $\mathcal{D}$ is essentially small and $\mathcal{C}$ is locally small and complete, then this Kan extension always exists. Furthermore, the Kleisli category of the codensity monad $\mathbb{T}$ can be seen as a full subcategory of $[\mathcal{D},\textbf{Set}]^\text{op}$~\cite{leinster_codensity_2013}, which gives us a \textit{formal setting} to interpret the Kleisli morphisms of the monad $\mathbb{T}$. We will give an account of the theory of codensity monads at the beginning of \hyperref[section3]{Section 3}. We then recall the presentation in~\cite{van_belle_probability_2022} of probability monads as codensity monads of functors from the category $\textbf{FinStoch}$ (or $\textbf{cStoch}$) of finite (respectively countable) sets and random functions.  Intuitively, a functor $K_P \colon \textbf{FinStoch} \to \mathcal{C}$ defines a \emph{model} of discrete probability in $\mathcal{C}$, and its codensity monad $P$ can be viewed as giving the largest extension of this notion. 
\subsection{Measure-theoretic probability from codensity}  

Measure-theoretic probability has a monadic formulation via the Giry monad, $\mathcal{G}$ on the category of measurable spaces, \textbf{Meas}~\cite{giry_categorical_1982}. If $P$ is a probability monad on a category $\mathcal{C}$ described by \hyperref[table1]{Table 1}, we expect a connection between $P$ and $\mathcal{G}$ which allows us to formulate $P$ using measure theory. In \hyperref[section4]{Section 4}, we develop this connection by showing that many probability monads are \emph{universal liftings} of the Giry monad. In particular, when $\mathcal{C}$ is a category of topological or metric spaces, there is a faithful functor $H \colon \mathcal{C} \to \textbf{Meas}$ which assigns a measurable structure to objects of $\mathcal{C}$. In this case, an endofunctor $P \colon \mathcal{C}\to \mathcal{C}$ with a natural transformation $HP \to \mathcal{G}H$ is an \emph{endofunctor of probability measures}, since we can interpret elements of $HPX$ as probability measures on $HX$ via this natural transformation.

Now suppose $P$ is given by the codensity monad of a functor $K_P \colon \textbf{FinStoch} \to \mathcal{C}$. Then $K_P$ defines a model of discrete probability in $\mathcal{C}$, $P$ defines its largest admissible extension, and $H$ defines its \textit{measurable interpretation}. Suppose $H$ and $K_P$ satisfy certain compatibility conditions, which ensure that the measurable interpretation of probability in $\mathcal{C}$ corresponds to measure-theoretic probability. In that case, we show in \hyperref[thm4.9]{Theorem 4.9} that codensity monads of $K_P$ have an additional universal property as the \textit{largest lifting} of the Giry monad. In particular, the natural transformation $HP \to \mathcal{G}H$ will additionally preserve the monad structure in this case. \hyperref[thm4.9]{Theorem 4.9} also provides the converse; it shows that universal liftings of a codensity monad inherit a codensity presentation.

\hyperref[thm4.9]{Theorem 4.9} is shown to generalise a result in \cite{van_breugel_metric_nodate}, which has a more limited universal lifting property for the Kantorovich monad $\mathcal{K}$ on the category of compact metric spaces, $\textbf{KMet}$. It also provides a novel characterisation of the Radon monad $\mathcal{R}$ on the category of compact Hausdorff spaces, $\textbf{KHaus}$. The compatibility conditions for the functors $H$ in these examples are non{-}trivial and rely on results proven in \hyperref[appenA]{Appendix A}. This appendix develops theory for the Baire functor $H \colon \textbf{KHaus} \to \textbf{Meas}$, which is shown in \hyperref[thmA.3]{Theorem A.3} to be the unique functor (up to isomorphism) that preserves limits and such that $H[0,1]\cong [0,1]$.

\subsection{Monoidal structure of probability monads}  

Categorical descriptions of probability theory rely on the additional structure of a \textit{lax monoidal monad}~\cite{kock_strong_1972} in order to model concepts such as marginals and independence. Explicitly, this requires a map $\chi_{X,Y} \colon PX\otimes PY\to P(X\otimes Y)$ that is coherent with the monad structure. In the case of the probability monads described in \hyperref[section2]{Section 2}, this is uniquely determined and sends a pair of measures $(p,q)$ to the \emph{product measure} $p\otimes q$. This can also be understood to be the unique map satisfying, for maps $h \colon X\to [0,1],k \colon Y\to [0,1]$,
\[
\int_{X\times Y} h \cdot k \; \text{d}(p\otimes q) = \int_X h \; \text{d}p\cdot \int_Y k \; \text{d}q
\]
This description of the monoidal structure in terms of integration operators motivates us to introduce \textit{exactly pointwise monoidal} codensity monads in \hyperref[section5]{Section 5}. These are codensity monads that satisfy additional conditions, appropriate for the monoidal setting. In particular,
\begin{enumerate}
    \item They have a lax monoidal structure, defined in terms of their limit presentation.
    \item The main result of this section, \hyperref[thm5.23]{Theorem 5.23}, shows that a codensity monad is exactly pointwise monoidal if and only if its Kleisli category can be seen as a monoidal subcategory of $[\mathcal{D},\textbf{Set}]^\text{op}$ in a suitable sense. 
\end{enumerate}
This theorem justifies our study of exactly pointwise monoidal codensity monads. It extends the \textit{formal setting} for interpreting Kleisli morphisms to monoidal categories, and gives a description of the tensor product of free algebras of these monads in terms of Day convolution.

In order to describe these additional conditions for probability monads, we develop some theory of probability $k${-}polymeasures. These are functions that take tuples of measurable subsets $(A_1,\dots,A_k)$, which are measures separately in each variable. We show that the exactly pointwise monoidal conditions require that the space of probability $k${-}polymeasures coincides with the space of measures on the product. In \hyperref[example5.15]{Example 5.15} we prove that this condition holds for the Radon monad. In \hyperref[example5.16]{Example 5.16} we show that this condition holds for the Giry monad when restricted to standard Borel spaces, but does not hold for arbitrary measurable spaces.

\subsection{Notation}  

\textbf{Set} denotes the category of sets in ZFC and functions. If $A \subseteq X$ is any subset, we denote its characteristic function by $\mathbbm{1}_A \colon X \to [0,1]$, but if it is clear from context, we may use $\mathbbm{1}_A$ to denote the characteristic function valued in $\{0,1\}$. \textbf{Meas} is the category of measurable spaces and measurable maps. If $X$ is a measurable space, $\mathcal{B}X$ is its $\sigma${-}algebra. If $X$ is a topological or metric space, $\mathcal{B}X$ is its Borel $\sigma${-}algebra. \textbf{KHaus} is the category of compact Hausdorff spaces and continuous maps. If $X,Y$ are topological spaces, $C(X,Y)$ is the set of continuous maps $X \to Y$ and $C(X)$ is the set of continuous maps into $\mathbb{R}$. \textbf{KMet} is the category of compact metric spaces and non{-}expansive maps, which are $f \colon X \to Y$ satisfying $\text{d}_Y(f(x),f(y))\leq  \text{d}_X(x,y)$. \textbf{Cat} denotes the 2{-}category of locally small categories and \textbf{MonCat} denotes the 2{-}category of locally small monoidal categories, lax monoidal functors, and monoidal natural transformations. If $X_i$ are objects of a monoidal category for $i=1,\dots,k$ then $\bigotimes_{i=1}^k X_i=(\dots(X_1\otimes_\mathcal{C} X_2)\otimes_\mathcal{C}\dots)\otimes_\mathcal{C}X_k$ is their $k${-}fold monoidal product associated to the left, and is defined similarly on morphisms. If $k=0$, then denotes the monoidal unit $I$. \textbf{SMonCat} denotes the locally full sub{-}2{-}category of symmetric monoidal categories and symmetric monoidal functors. If a definition references a monoidal structure and the underlying category has finite products, then it is given by the Cartesian structure unless otherwise stated.

\subsection{Acknowledgements}  

I am grateful to the many people who have given helpful feedback. In particular, I would like to thank Sam Staton and Bartek Klin for many helpful discussions. I would also like to thank all of the people in Sam Staton's group with whom I've discussed this, with particular thanks to Ruben Van Belle, who talked through the analytical issues with me and suggested the restriction to \textbf{BorelMeas} in 
 \hyperref[example5.16]{Example 5.16}, and Paolo Perrone, who gave useful pointers and clarifications. I'm also grateful to the anonymous reviewer, who gave many helpful suggestions.
\section{Five probability monads}
\label{section2}
 In this section, we define the five main examples of probability monads in this paper.
\subsection{Countable distribution monad $\mathcal{D}_c$ on \textbf{Set}}  

  Define $\mathcal{D}_c \colon \textbf{Set} \to \textbf{Set}$ by \[\mathcal{D}_cX=\{p \colon X \to [0,1] : \sum_{x \in X} p(x)=1\}\] and for $f \colon Y \to X$ by \[\mathcal{D}_cf(p)(x)=\sum_{y \in f^{{-}1}(\{x\})}p(y)\] Note that this definition implies that $p\in \mathcal{D}_cX$ must have countable support. Now, $\eta_X \colon X \to \mathcal{D}_cX$, defined by \[ \eta_X(x)(y)=\mathbbm{1}_{\{x\}}(y)\] and $\mu_X \colon \mathcal{D}_c \mathcal{D}_c X \to \mathcal{D}_cX$, defined by \[\mu_X(\mathbb{P})(x)=\sum_{p \in \mathcal{D}_cX} \mathbb{P}(p)\cdot p(x)\] provide the components of natural transformations that give $\mathcal{D}_c$ a monad structure.

  $\mathcal{D}_c$ has a lax monoidal structure  $\chi_{X,Y} \colon \mathcal{D}_cX \times \mathcal{D}_cY \to \mathcal{D}_c(X \times Y)$, defined by \[\chi_{X,Y}(p,q)(x,y)=p(x)\cdot q(y)\] $\mathcal{D}_c$ is \emph{affine}.
    \begin{definition}[Affine monad]
    A monad $\mathbb{T}=(T,\eta,\mu)$ on a category $\mathcal{C}$ with a terminal object 1 is affine if $T1 \cong 1$
  \end{definition}
  
  A finitary variant of this monad is studied in \cite{jacobs_probability_2018}, and its algebras are discussed in \cite{fritz_convex_nodate} and originate in \cite{marshall_h_stone_postulates_1949}.
\subsection{Giry monad $\mathcal{G}$ on \textbf{Meas}~\cite{giry_categorical_1982}}  \label{example2.3}
 
   Define $\mathcal{G} \colon \textbf{Meas} \to \textbf{Meas} $ by \[\mathcal{G}X=\{p \colon \mathcal{B}X \to [0,1] : p \;\text{is a probability measure}\}\] with the coarsest $\sigma${-}algebra such that $\text{ev}_A \colon \mathcal{G}X \to [0,1]$, where $\text{ev}_A(p)=p(A)$, is measurable for each $A \in \mathcal{B}X$. Note that since every measurable function is a monotone limit of simple functions, for any measurable $h \colon X \to [0,1]$, the function $\text{ev}_h \colon \mathcal{G}X \to [0,1]$ sending $p\mapsto \int_X h \; \text{d}p$ is also measurable. Hence, the $\sigma${-}algebra on $\mathcal{G}X$ is also initial for the maps $\text{ev}_h$. For $f \colon Y \to X$, define $\mathcal{G}f$ to send $p$ to its pushforward measure along $f$, defined by \[\mathcal{G}f(p)(A)=p(f^{{-}1}(A))\] This is measurable since $\text{ev}_h \mathcal{G}f=\text{ev}_{hf}$. Define $\eta_X \colon X \to \mathcal{G}X$ to send $x$ to its Dirac measure, defined by \[\eta_X(x)(A)=\mathbbm{1}_A(x)\] and define $\mu_X \colon \mathcal{G}\mathcal{G}X \to \mathcal{G}X$ to send $\mathbb{P}$ to its average value, defined by \[\mu_X(\mathbb{P})(A)=\int_{\mathcal{G}X}\text{ev}_A \; \text{d}\mathbb{P}\] These are measurable since $\text{ev}_h \eta_X=h$ and $\text{ev}_h\mu_X=\text{ev}_{\text{ev}_h}$, and they give $\mathcal{G}$ the structure of a monad~\cite{giry_categorical_1982}.

  $\mathcal{G}$ has a lax monoidal structure $\chi_{X,Y} \colon \mathcal{G}X \times \mathcal{G}Y \to \mathcal{G}(X \times Y)$ given by $\chi_{X,Y}(p,q)=p\otimes q$, the unique measure satisfying \[\text{ev}_{A\times B}(p\otimes q)=\text{ev}_A(p)\cdot\text{ev}_B(q)\] The commutativity of $\mathcal{G}$ then follows from Fubini's theorem for product measures. $\mathcal{G}$ is also affine.

\subsection{Countable expectation monad $\mathcal{E}_c$ on \textbf{Set}} \label{example2.4}

  Define $\mathcal{E}_c:\textbf{Set} \to \textbf{Set}$ by $\mathcal{E}_c=U\mathcal{G}D$, where $D \colon \textbf{Set}\to \textbf{Meas}$ assigns a set its discrete $\sigma${-}algebra and $U \colon \textbf{Meas}\to \textbf{Set}$ is the forgetful functor. Explicitly, \[\mathcal{E}_cX=\{p:\mathcal{P}X \to [0,1] : p \;\text{is a probability measure}\}\] where $\mathcal{P}X$ is the power set of $X$. In fact, it is consistent with ZFC that $\mathcal{E}_c=\mathcal{D}_c$~\cite{ulam_zur_1930} since there are models of ZFC where the only measures on $\mathcal{P}X$ are given by discrete distributions~\cite[Chapter 54]{fremlin_measure_2000}. However, in other models of ZFC, $\mathcal{E}_c$ may not be commutative.
  \begin{proposition}
      The commutativity of $\mathcal{E}_c$ is undecidable in \textnormal{ZFC}.
  \end{proposition}
  \proof
  If $\mathcal{E}_c=\mathcal{D}_c$, it is clearly commutative. However, if $\mathcal{E}_c \neq \mathcal{D}_c$ then there exists a smallest measurable cardinal $\kappa$~\cite{fremlin_measure_2000}. Let $\beta_\kappa$ be the submonad of $\mathcal{E}_c$ of 0{-}1 valued measures, explicitly \[\beta_\kappa X = \{p\in \mathcal{E}_cX :   \forall A \in \mathcal{P}X \; p(A) \in \{0,1\} \}\] In fact $\beta_\kappa$ is isomorphic to the submonad of the ultrafilter monad~\cite{leinster_codensity_2013,kennison_equational_1971,manes_algebraic_1976} of $\kappa${-}complete ultrafilters. Now, if $\mathcal{U}$ is a non{-}trivial $\kappa${-}complete ultrafilter on $\kappa$ then the left and right strengths of $\beta_\kappa$ induce two $\kappa${-}complete ultrafilters on $\kappa \times \kappa$ given by \[\mathcal{U}\otimes_r \mathcal{U}=\{M \subseteq \kappa \times \kappa \colon \{x \in \kappa \colon \{y \in \kappa \colon (x,y) \in M \} \in \mathcal{U} \} \in \mathcal{U} \} \] and similarly for $\mathcal{U}\otimes_l \mathcal{U}$. But the set $\Delta^{-} =\{(x,y) \in \kappa \times \kappa \colon x<y \}$ satisfies $\Delta^{-} \in \mathcal{U}\otimes_r \mathcal{U}$ but $\Delta^{-} \notin \mathcal{U}\otimes_l \mathcal{U}$. Hence, $\beta_\kappa$ is not commutative and thus neither is $\mathcal{E}_c$.
  \endproof
  Intuitively, this is because the $\sigma${-}algebra on the product of uncountable discrete spaces is not discrete. If we let $\mathcal{P}X \otimes \mathcal{P}Y$ be the $\sigma$-algebra generated by all rectangles, then in general $\mathcal{P}X \otimes \mathcal{P}Y \neq \mathcal{P}(X\times Y)$. A variant of this monad, which takes finitely additive measures on the power set, was described in \cite{jacobs_probability_2018}.
\subsection{Radon monad $\mathcal{R}$ on \textbf{KHaus}~\cite{semadeni_monads_1973,tadeusz_swirszcz_monadic_1974}.}
  For $X$ a compact Hausdorff space, define \[\mathcal{R}X = \{ p \colon \mathcal{B}X \to [0,1] : p \text{ is a Radon probability measure}\}\]  recalling that a Borel probability measure on a compact Hausdorff space is Radon~\cite{fedorchuk_probability_1991} if \[p(A)=\sup\{p(K) : K \subseteq A \; \text{is compact}\}\]  for any Borel set $A$. Define the topology on $\mathcal{R}X$ to be the coarsest such that $\text{ev}_h$ is continuous for any continuous $h \colon X \to [0,1]$. 
  \begin{proposition}[\cite{semadeni_monads_1973,tadeusz_swirszcz_monadic_1974}]
   $\mathcal{R}X$ is a compact Hausdorff space.
  \end{proposition}
  \proof
  Let $C(X)^*$ be the space of linear functionals on $X$ with the weak-$*$ topology. By the Riesz{-}Markov representation theorem~\cite{rudin_functional_1991}, $\mathcal{R}X$ is homeomorphic to \[M^+(X)\cap B_1 \cap \{\mu : \text{ev}_X(\mu)=1\}\]
  Here, $B_1$ is the closed unit ball of $C(X)^*$, and $M^+(X)$ is the cone of positive measures, given by
  \[M^+(X)=\bigcap_{f\in \mathcal{C}(X), f\geq 0}\{\mu :\text{ev}_f(\mu) \geq 0\}\]
  In particular, $M^+(X)$ is closed, and so $\mathcal{R}X$ is homeomorphic to a closed subspace of the unit ball $B_1$. By the Banach-Alaoglu theorem, the unit ball is a compact Hausdorff space~\cite{rudin_functional_1991}, and so $\mathcal{R}X$ is also a compact Hausdorff space.
\endproof
We then define $\mathcal{R}$ on morphisms via the pushforward measure, and its multiplication and unit analogously to \hyperref[example2.3]{Section 2.3}. That these maps are continuous is analogous to \hyperref[example2.3]{Section 2.3}. The pushforward of a Radon measure is also Radon, Dirac measures are Radon, and the multiplication preserves Radon measures. Hence, $\mathcal{R}$ has the structure of a monad on \textbf{KHaus}. The algebras of the Radon monad are compact and convex subspaces of locally convex topological vector spaces~\cite{keimel_monad_2008}.
\subsection{Kantorovich monad $\mathcal{K}$ on \textbf{KMet}}  

  For $X$ a compact metric space, we define \[\mathcal{K}X= \{p \colon \mathcal{B}X \to [0,1] : p  \; \text{is a probability measure}\}\] and \[\text{d}_{\mathcal{K}X}(p,q)=\sup\{\left| \int_X h \; \text{d}p -\int_X h \; \text{d}q \right| : h \colon X \to [0,1] \text{ is non{-}expansive} \}\] This is the Kantorovich metric (also known as the 1{-}Wasserstein distance) used in optimal transport theory~\cite{basso_hitchhikers_2015}. Note that every Borel probability measure on a compact metric space is Radon (see \hyperref[appenA]{Appendix A}). Also, the topological space induced by $\text{d}_{\mathcal{K}X}$ on $\mathcal{K}X$ coincides with $\mathcal{R}X$, when $X$ is assigned its metric topology. Now, $\text{ev}_h$ is non{-}expansive for each non{-}expansive $h \colon X \to [0,1]$, and a map $f \colon Y \to \mathcal{K}X$ is non{-}expansive iff $\text{ev}_h f$ is non{-}expansive for each non{-}expansive $h$. Thus, we define $\mathcal{K}f$ for $f \colon Y \to X$, $\eta$ and $\mu$ as in previous examples, and $\mathcal{K}$ has the structure of a monad on \textbf{KMet}. There is a symmetric monoidal structure on \textbf{KMet} where $(X,d_X)\otimes (Y,d_Y)=(X\times Y,d_{X \otimes Y})$ and $d_{X \otimes Y}((x,y),(x',y'))=d_X(x,x')+d_Y(y,y')$. Note that this monoidal structure is not the cartesian monoidal structure. With this structure, the Kantorovich monad is shown to have a lax monoidal structure in \cite{fritz_bimonoidal_2018}. This monad was introduced in \cite{van_breugel_metric_nodate} and is further studied in \cite{fritz_probability_2019}.
\begin{remark}
It is possible to extend some of these monads to monads of finitely additive measures, subprobability measures, and other generalisations. However, this can be at the cost of affineness and commutativity. In this paper, we restrict our view to the monads of probability and subprobability measures, which we will study in \hyperref[prop3.12]{Proposition 3.12}.
\end{remark}
\section{Codensity monads and probability} \label{section3}
In this section, we give an overview of the theory of codensity monads, and show how the probability monads from \hyperref[section2]{Section 2} can be presented as pointwise codensity monads of functors from essentially small categories of stochastic maps. In the final part of this section, we give sufficient conditions for a codensity monad to be affine, which are satisfied by the monads in \hyperref[section2]{Section 2}. We then show that, under the same conditions, monads of subprobability measures inherit limit descriptions from their probability counterparts.
\subsection{Codensity monads}
\begin{definition}[Codensity monad] \label{definition3.2} Let $K \colon \mathcal{D} \to \mathcal{C}$ be any functor such that there is a diagram
\[\begin{tikzcd}
	{\mathcal{D}} & {\mathcal{C}} \\
	{\mathcal{C}} 
	\arrow["K", from=1-1, to=1-2]
	\arrow[""{name=0, anchor=center, inner sep=0}, "K"', from=1-1, to=2-1]
	\arrow[""{name=1, anchor=center, inner sep=0}, "T", from=1-2, to=2-1]
	\arrow["\varepsilon"', yshift = 1mm, shorten <=3pt, shorten >=3pt, Rightarrow, from=1, to=0]
\end{tikzcd}\]
 with $(T,\varepsilon)$ as the right Kan extension of $K$ along itself~\cite{mac_lane_categories_1998}. Then $T$ alongside the unique maps $\mu$ and $\eta$ satisfying $\varepsilon\, \mu_K=\varepsilon \, T\varepsilon$ and $\varepsilon \, \eta_K=1_K$ define a monad structure $\mathbb{T}$, and $(\mathbb{T},\varepsilon)$ is called the \textit{codensity monad of $K$}. In the case that the Kan extension is pointwise, we call $(\mathbb{T},\varepsilon)$ a \textit{pointwise codensity monad}.
 \end{definition}
  Pointwise codensity monads inherit a limit presentation from the limit formula for pointwise Kan extensions~\cite{mac_lane_categories_1998}. Namely, \[TA \cong \lim \, (A\downarrow K) \xrightarrow{U_A} \mathcal{D} \xrightarrow{K} \mathcal{C}\] where $(A\downarrow K)$ is a comma category. Written as an end, $TA \cong \int_{\mathcal{D}}K{-}^{\mathcal{C}(A,K{-})}$, where if $A \in \text{ob} \, \mathcal{C}$ and $X$ is a set, $A^X$ denotes the $X$th power of $A$ in $\mathcal{C}$. In particular, if $\mathcal{C}=\textbf{Set}$ and $K \colon \mathcal{D} \to \textbf{Set}$ we have $ TA \cong [\mathcal{D},\textbf{Set}](\mathcal{C}(A,K{-}),K) $, so $\mathbb{T}$ is, in a sense, a \textit{double dualisation monad}. We also obtain an explicit description of the unit and multiplication from this limit presentation. Since $TA$ is a limit over the comma category $(A \downarrow K)$, for any $h \colon A \to KX$, there is a cone leg $\text{ev}_h \colon TA \to KX$, where explicitly $\text{ev}_h = \varepsilon_X Th$. Now, for $f \colon B \to A$, $Tf$ is the unique map satisfying $\text{ev}_h Tf = \text{ev}_{hf}$,  $\eta_A \colon A \to TA$ is the unique map satisfying $\text{ev}_h \eta_A = h$, and $\mu_A \colon TTA \to TA$ is the unique map satisfying $\text{ev}_h \mu_A = \text{ev}_{\text{ev}_h}$. In the following setting, codensity monads are always guaranteed to exist and be pointwise; hence the distinction between pointwise and non{-}pointwise codensity monads is dropped under these assumptions:
 
 \begin{proposition}[\cite{leinster_codensity_2013}] \label{prop3.3}
      If $\mathcal{D}$ is essentially small, $\mathcal{C}$ is locally small and complete, and $K \colon \mathcal{D} \to \mathcal{C}$ is a functor, then there is a right adjoint to the functor $K_{\circ} \colon \mathcal{C} \to [\mathcal{D},\textnormal{\textbf{Set}}]^\textnormal{\text{op}} $ given by currying $\mathcal{C}({-},K{=})$. This right adjoint factorises $K$ through $[\mathcal{D},\textnormal{\textbf{Set}}]^\textnormal{\text{op}}$.
      \[\begin{tikzcd}
	{\mathcal{D}} & {[\mathcal{D},\textnormal{\textbf{Set}}]^\textnormal{\text{op}}} \\
	& {\mathcal{C}}
	\arrow["{Y_\mathcal{D}}", from=1-1, to=1-2]
	\arrow["K"', from=1-1, to=2-2]
	\arrow[""{name=0, anchor=center, inner sep=0}, "{\textnormal{Ran}_{Y_\mathcal{D}}K}", shift left=2, from=1-2, to=2-2]
	\arrow[""{name=1, anchor=center, inner sep=0}, "{K_\circ}", shift left=2, from=2-2, to=1-2]
	\arrow["\dashv"{anchor=center}, draw=none, from=1, to=0]
\end{tikzcd}\]
Furthermore, the monad induced by this adjunction is the pointwise codensity monad of $K$.
 \end{proposition} 
In particular, the Kleisli comparison functor gives a full and faithful embedding $\mathcal{C}_\mathbb{T} \to [\mathcal{D},\textbf{Set}]^\text{op}$. Hence, we can exhibit $\mathcal{C}_\mathbb{T}$ as the full subcategory of functors of the form $\mathcal{C}(A,K{-})$ for some $A\in \text{ob} \; \mathcal{C}$.

It is also worth noting that every monad arises as a pointwise codensity monad of some functor. If $G$ has a left adjoint $F$, then the codensity monad of $G$ is given by the monad structure on $GF$ induced by the adjunction. In this case, we have $\text{Ran}_GG \cong G \, \text{Ran}_G 1$, and since the latter is an absolute Kan extension, the former is pointwise, so every monad arises as a pointwise codensity monad in this way. For this reason, we do not consider codensity monads as a particular class of monads, or even a construction of a monad, but rather as a convenient presentation of a monad. When $K \colon \mathcal{D} \to \mathcal{C}$ is a functor from an essentially small category to a complete locally small category, one can present $\mathbb{T}$ in terms of simple data about the essentially small categories $\mathcal{D}$ and $(A \downarrow K)$ rather than complex data about the categories $\mathcal{C}$ or $[\mathcal{D},\textbf{Set}]^{\text{op}}$.
\subsection{Codensity presentations of probability monads}  

 One sensible aim, given a monad $\mathbb{T}$ on $\mathcal{C}$, might be to find an essentially small category $\mathcal{D}$ and a functor $K \colon \mathcal{D} \to \mathcal{C}$ such that $\mathbb{T}$ is the pointwise codensity monad of $K$. This was done for the probability monads in \hyperref[section2]{Section 2} in \cite{van_belle_probability_2022}. Given this data, we could then try to express additional structure and properties of $\mathbb{T}$ in terms of its limit presentation, and the structure and properties of $K$. It is sufficient to look at subcategories of the category of Eilenberg{-}Moore algebras $\mathcal{C}^\mathbb{T}$, since any such $K$ must factor through the Eilenberg{-}Moore forgetful functor $G^\mathbb{T} \colon \mathcal{C}^\mathbb{T} \to \mathcal{C}$~\cite{leinster_codensity_2013}. Furthermore, if $\mathbb{T}$ is the pointwise codensity monad of a functor $K \colon \mathcal{D} \to \mathcal{C}$, it is also the pointwise codensity monad of the inclusion of the smallest subcategory containing all the morphisms $Km$ for $m \in \text{mor} \; \mathcal{D}$, and so one may consider faithful $K$. 
  Since $\mathbb{T}$ is also the codensity monad of the Kleisli forgetful functor $G_\mathbb{T}$, it is reasonable to look at the restrictions of $G_\mathbb{T}$ to subcategories of the Kleisli category $\mathcal{C}_\mathbb{T}$. 
  \begin{definition}[Kleisli codensity presentation]
      Let $I \colon \mathcal{D} \to \mathcal{C}_\mathbb{T}$ be the inclusion of a subcategory. $\mathbb{T}$ has a Kleisli codensity presentation over $\mathcal{D}$ if, for each  $A\in \textnormal{ob} \; \mathcal{C}$, $(TA,\mu_{TX}Th)$ is a limit cone over $(A\downarrow G_\mathbb{T}I)$.
  \end{definition}

If $\mathbb{T}$ has a Kleisli codensity presentation over $\mathcal{D}$, then it additionally has a Kleisli codensity presentation over the full subcategory of $\mathcal{C}_\mathbb{T}$ on objects in $\mathcal{D}$. We now provide Kleisli codensity presentations, over essentially small categories, for the probability monads we described in \hyperref[section2]{Section 2}.
\begin{definition} The category \textnormal{\textbf{FinStoch}} has finite sets as its objects and random functions as its morphisms~\cite{fritz_synthetic_2020}. It is the full subcategory of the Kleisli category $\textnormal{\textbf{Set}}_{\mathcal{D}_c}$ on finite sets. \textnormal{\textbf{cStoch}} is the full subcategory of $\textnormal{\textbf{Set}}_{\mathcal{D}_c}$ on countable sets.
\end{definition}
 Since $\mathcal{D}_c$ is a lax monoidal monad, both categories have a symmetric monoidal structure which is given on objects by the Cartesian product of sets. Additionally, there are faithful, identity on object functors as in the following diagram
\[\begin{tikzcd}
	{\textbf{FinSet}} & {\textbf{FinStoch}} \\
	{\textbf{cSet}} & {\textbf{cStoch}}
	\arrow["{{F_{\mathcal{D}_c,f}}}", from=1-1, to=1-2]
	\arrow[hook, from=1-1, to=2-1]
	\arrow[hook, from=1-2, to=2-2]
	\arrow["{{F_{\mathcal{D}_c,c}}}"', from=2-1, to=2-2]
\end{tikzcd}\]
given by the restriction of the Kleisli left adjoint $F_{\mathcal{D}_c}$, where $\textbf{cSet}$ is the category of countable sets.

\begin{example}
\textbf{Giry monad on Meas.} \label{example3.7} 

There is a full and faithful functor $ \textbf{cStoch}\to \textbf{Meas}_\mathcal{G}$ defined by assigning a countable set its discrete $\sigma${-}algebra. Let $K_\mathcal{G}$ be the restriction of the Kleisli forgetful functor $G_\mathcal{G}$ to $\textbf{cStoch}$. It is shown in \cite{van_belle_probability_2022} that $\mathcal{G}$ is the pointwise codensity monad of the composite functor $K_\mathcal{G} F_{\mathcal{D}_c,c} \colon \textbf{cSet}\to \textbf{Meas}$. Since $F_{\mathcal{D}_c,c}$ is a faithful, identity on objects functor, it also follows that $\mathcal{G}$ is the pointwise codensity monad of the functor $K_\mathcal{G}$, and has a Kleisli codensity presentation over \textbf{cStoch}.  We now sketch the proof that the Giry monad is the codensity monad of $K_\mathcal{G} F_{\mathcal{D}_c,c}$~\cite{van_belle_probability_2022}.

The cone $(\mathcal{G}A, \text{ev}_h)$ over $(A\downarrow K_\mathcal{G} F_{\mathcal{D}_c,c})$ with legs $\text{ev}_h=\mu_X\mathcal{G}h$, for $h \colon A \to K_\mathcal{G} F_{\mathcal{D}_c,c}X$, is explicitly defined by  $\text{ev}_h(p)(x)=\int h({-})(x) \; \text{d}p$. Here we implicitly use the isomorphism \[K_\mathcal{G} F_{\mathcal{D}_c,c}X \cong\{p \in [0,1]^X: \sum_{x\in X}p(x)=1 \}\] with $\sigma${-}algebra given as a subspace of the product of the Borel structure on $[0,1]$. Now if $(C,\lambda_h)$ is any other cone, then for any $B\in \mathcal{B}A$ we can define $h_B= (\mathbbm{1}_B,\mathbbm{1}_{B^c}) \colon A \to K_\mathcal{G} F_{\mathcal{D}_c,c}\{0,1\}$ and $\lambda_B=\lambda_{h_B}({-})(0)\colon C \to [0,1]$. Then, for any $c \in C$, $\alpha(c)(B)=\lambda_B(c)$ defines a probability measure on $A$, giving a unique factorisation of the cone maps, and demonstrating that $(\mathcal{G}A, \text{ev}_h)$ is a limit cone. 
\end{example}

\begin{example}
\textbf{Countable expectation monad on Set.}

First, we note that the restriction of $\mathcal{E}_c$ to $\textbf{cSet}$ is isomorphic to the restriction of $\mathcal{D}_c$ to \textbf{cSet}, since all countably supported measures are discrete. Thus, the subcategory of $\textbf{Set}_{\mathcal{E}_c}$ with objects in $\textbf{cSet}$ is (isomorphic to) $\textbf{cStoch}$. We write $K_{\mathcal{E}_c} \colon \textbf{cStoch} \to \textbf{Set}$ for the restriction of the Kleisli forgetful functor. Then $\mathcal{E}_c$ has a Kleisli codensity presentation over $\textbf{cStoch}$, since it is the pointwise codensity monad of $K_{\mathcal{E}_c}$. This can be shown by a similar proof as used in \hyperref[example3.7]{Example 3.7}. However, we will prove this in \hyperref[exmp4.16]{Example 4.16} as a corollary of \hyperref[thm4.9]{Theorem 4.9}.
\end{example}

\begin{example}
\label{example3.9}
\textbf{Radon monad on KHaus, Kantorovich monad on KMet.}

There is a full and faithful functor $ \textbf{FinStoch} \to \textbf{KHaus}_\mathcal{R}$, assigning sets their discrete topology. We write $K_\mathcal{R} \colon \textbf{FinStoch} \to \textbf{KHaus}$ for the restriction of the Kleisli forgetful functor to this subcategory. Then, it is shown in \cite{van_belle_probability_2022} that the Radon monad is the codensity monad of $K_\mathcal{R} F_{\mathcal{D}_c,f}$. Hence, the Radon monad has a Kleisli codensity presentation over $\textbf{FinStoch}$, since it is also the pointwise monad of $K_\mathcal{R}$. Similarly, by analogous results in \cite{van_belle_probability_2022}, the Kantorovich monad has a Kleisli codensity presentation over $\textbf{FinStoch}$, and we write the restriction of the Kleisli forgetful functor as $K_\mathcal{K}$ in this case. 
\end{example}
\subsection{Subprobability measures and affineness}  

 All the probability monads we have encountered are affine. There is only one probability measure on a singleton, since the total probability of events is 1. One can also use the codensity presentation of these monads to show this.
\begin{proposition}
\label{prop3.11}
    Suppose $\mathcal{C}$ and $\mathcal{D}$ have terminal objects. If $K \colon \mathcal{D} \to \mathcal{C}$ is such that $!_{K1} \colon K1 \to 1$ is an isomorphism, and for every $y \colon 1 \to KX$ there is a unique $x \colon 1 \to X$ with $y=Kx \, !_{K1}^{{-}1}$. Then, if $T$ is the pointwise codensity monad of $K$, $\eta_1$ is an isomorphism. 
\end{proposition}
\proof
    Since $1$ is terminal, $!_{T1}\eta_1=1_1$, so it is sufficient to show $\eta_1!_{T1}=1_{T1}$. For each $y \colon 1 \to KX$, $\text{ev}_y \eta_1 !_{T1}=y!_{T1}$ but $y=Kx!_{K1}^{{-}1}$ by assumption, and $!_{K1}^{{-}1}!_{T1}=\text{ev}_{!_{K1}^{{-}1}} \colon T1 \to K1$ so $y!_{T1}=Kx \text{ev}_{!_{K1}^{{-}1}}=\text{ev}_y$. Hence, $\eta_1!_{T1}=1_{T1}$ by the uniqueness of cone factorisation.
\endproof

 For locally small categories $\mathcal{D}$ and $\mathcal{C}$, the latter condition is equivalent to the map $\mathcal{D}(1,{-}) \to \mathcal{C}(1,K{-})$ being an isomorphism, but the argument above only relied on it being epic. If $\mathbb{T}$ is an affine monad, and $\mathcal{D}$ is a full subcategory of $\mathcal{C}_\mathbb{T}$ containing 1, then the restriction of $G_\mathbb{T}$ to $\mathcal{D}$ satisfies this condition. In particular, the functors $K_\mathcal{G}$, $K_{\mathcal{E}_c}$, $K_\mathcal{R}$ and $K_\mathcal{K}$ satisfy the conditions in \hyperref[prop3.11]{Proposition 3.11}, noting that the terminal object in $\textbf{FinStoch}$ and $\textbf{cStoch}$ is the singleton set 1 since $D$ is affine. However, applying this argument to an individual monad is somewhat circular; if $\mathcal{D}$ is a full subcategory of $\mathcal{C}_\mathbb{T}$ with a terminal object not of the form $F_\mathbb{T}1$, the latter condition may not hold (consider, for example, the power set monad $\mathcal{P}$ on \textbf{Set}).
 
If $\mathcal{C}$ has a terminal object $1$ such that $1+X$ exists for each $X$, then for any monad $\mathbb{T}$ on $\mathcal{C}$ there is a unique distributive law $\lambda \colon 1+T{-} \to T(1+{-})$~\cite{barr_toposes_1985} which induces a monad structure on the endofunctor $T(1+{-})$ ~\cite{van_breugel_metric_nodate}. In the case of the probability monads from \hyperref[section2]{Section 2}, we call this the \textit{monad of subprobability measures}. If  $p \in T(1+X)$ then $p(X)=1-p(\bot)$ and so elements of $T(1+X)$ correspond to measures on $X$ satisfying $p(X)\leq 1$. In a programming language setting, we can view this monad as combining probabilistic nondeterminism with possible failure, and $p(\bot)$ can be interpreted as the probability of an error. By allowing the cone maps $\text{ev}_h$ to incorporate possible failure, we can obtain a pointwise Kan extension expressing monads of subprobability measures:

\begin{proposition}
\label{prop3.12}
Under the conditions of \hyperref[prop3.11]{Proposition 3.11}, if binary coproducts with $1$ exists in $\mathcal{C}$ and $\mathcal{D}$, then $T(1+{-})$ is the pointwise right Kan extension of $K(1+{-})$ along $K$
\end{proposition}
\proof
First, we define the limit cones. If $h \colon A \to KX$, then $[K\iota_1 !_{K1}^{{-}1},K\iota_2 h] \colon 1 + A \to K(1+X)$ and for any $f \colon X \to Y$ we have \[K(1+f)\text{ev}_{[K\iota_1 !_{K1}^{{-}1},K\iota_2 h]}=\text{ev}_{[K\iota_1 !_{K1}^{{-}1},K\iota_2 Kf h]}\] Hence the maps $\text{ev}_{[K\iota_1 !_{K1}^{{-}1},K\iota_2 h]}\colon T(1+A) \to K(1+X)$ provide a cone over $K(1+U_A{-}) \colon (A \downarrow K) \to \mathcal{C}$ with apex $T(1+A)$.

Now if $(C,\tau_h)$ is any cone over $K(1+U_A{-})$ and $[y,h] \colon 1+A \to KX$, there is a unique map $x \colon 1 \to X$ such that $y=Kx!_{K1}^{{-}1}$ by hypothesis. Hence, we can define $(C,\nu_{[y,h]})$ by $\nu_{[y,h]}=K[x,1_X]\tau_h$ and for any $f \colon X \to Y$
\begin{align*}
Kf\nu_{[y,h]} & = K[fx,f]\tau_h \\
& =K[fx,1_X]\tau_{Kf h} \tag{$\tau$ is a cone}\\
& = \nu_{Kf[y,h]} \tag{$K(fx)!_{K1}^{{-}1}=Kfy$}
\end{align*}
Hence $(C,\nu_{[y,h]})$ is a cone over $KU_{1+A} \colon (1+A \downarrow K) \to \mathcal{C}$, and thus there is a unique map $\alpha \colon C \to T(1+A)$ such that $\text{ev}_{[y,h]}\alpha=\nu_{[y,h]}$. But $\alpha$ is also a map of cones over $(A \downarrow K)$, since
\begin{align*}
    \text{ev}_{[K\iota_1 !_{K1}^{{-}1},K\iota_2 h]}\alpha & =\nu_{[K\iota_1 !_{K1}^{{-}1},K\iota_2 h]} \\
    & = K[\iota_1,1]\tau_{K\iota_2 h} \\
    & = K[\iota_1,1]K(1+\iota_2)\tau_h \\
    & = \tau_h
\end{align*}
Finally, if $\alpha'$ is another cone map over $(A \downarrow K)$ such that $\text{ev}_{[K\iota_1 !_{K1}^{{-}1},K\iota_2 h]}\alpha'=\tau_h$, then
\begin{align*} \nu_{[y,h]} & =K[x,1_X]\text{ev}_{[K\iota_1 !_{K1}^{{-}1},F\iota_2 h]}\alpha' \\ & =\text{ev}_{[Kx!_{K1}^{{-}1},h]}\alpha'\\ & =\text{ev}_{[y,h]}\alpha' 
\end{align*}
hence $\alpha=\alpha'$, and so the cone factorisation is unique.
\endproof
\section{Probability monads as universal liftings and Kleisli laws} \label{section4}
 In this section, we discuss an alternative characterisation of the probability monads in \hyperref[section2]{Section 2} and relate it to their codensity presentations. We begin by giving an informal motivation. We can think of a terminal object in a category as the \textit{largest object satisfying some conditions} which arise from the uniqueness of the morphisms onto it. This gives an interpretation of probability monads expressed as the codensity monad of $K \colon \textbf{cStoch} \to \mathcal{C}$ as the \textit{largest endofunctor on $\mathcal{C}$ extending the model of discrete probability given by $K$}.
 
On the other hand, if an object of $\mathcal{C}$ can be assigned a measurable structure, we might expect a probability monad $\mathbb{T}$ on $\mathcal{C}$ to give a space of probability measures with respect to this structure. Explicitly, if $H \colon \mathcal{C} \to \textbf{Meas}$ is a functor, we expect a map $\alpha \colon HT \to \mathcal{G}H$ which assigns a measure on $HX$ for each element of $HTX$. If $(S,\beta \colon HS \to \mathcal{G}H)$ is any endofunctor of this form, elements of $HSA$ are probability measures on $HA$, and so $(S,\beta)$ can be viewed as an \emph{endofunctor of probability measures}. If $H$ and $K$ define compatible models of probability for $\mathcal{C}$, there are two informal perspectives on a maximal probability monad on $\mathcal{C}$ which we might expect to coincide: 
\begin{enumerate}
    \item $\mathbb{T}$ is the \textit{largest endofunctor extending discrete probability given by $K$} and so should be the codensity monad of $K$
    \item $\mathbb{T}$ is the \textit{largest endofunctor of probability measures with respect to $H$} and should be terminal in the category of endofunctors of probability measures.
\end{enumerate}
We make the connection between these notions precise in this section. In \hyperref[thm4.9]{Theorem 4.9}, we show that if certain compatibility conditions between $H$ and $K$ hold, then these two notions of maximal probability monads coincide. 
\subsection{Kleisli laws and monad liftings}
In order to relate monads on different categories, we recall the definition of a Kleisli law. 
\begin{definition}[Kleisli Law~\cite{mulry_lifting_1993}] \label{definition4.2} If $H \colon \mathcal{D} \to \mathcal{C}$ and $\mathbb{P}$ (resp. $\mathbb{T}$) is a monad on $\mathcal{D}$ (resp. $\mathcal{C}$), a \textit{Kleisli law} $\lambda \colon HP \to TH$ is a natural transformation satisfying
\begin{enumerate}
\item $\eta^\mathbb{T}_H=\lambda H\eta^\mathbb{P}$
\item $\mu^\mathbb{T}_H T\lambda \lambda_P=\lambda H\mu^\mathbb{P}$
\end{enumerate}
If $\mathcal{C}=\mathcal{D}$ and $H=1_\mathcal{C}$, then a Kleisli law $P \to T$ is referred to as a \textit{monad morphism}.
\end{definition}
\begin{proposition}[\cite{mulry_lifting_1993}]
\label{prop4.3}
For $H$, $\mathbb{P},\mathbb{T}$ as in \hyperref[definition4.2]{Definition 4.2}, there is a bijection between Kleisli laws $\lambda \colon HP \to TH$ and functors $\Bar{H} \colon \mathcal{D}_\mathbb{P} \to \mathcal{C}_\mathbb{T}$ satisfying $\Bar{H}F_\mathbb{P}=F_\mathbb{T}H$. Furthermore, such a functor $\Bar{H}$ satisfies $\Bar{H}f=\lambda_B Hf$ for $f \colon A \to B$ in $\mathcal{D}_\mathbb{T}$.
\end{proposition}
 Now if $\mathbb{T}$ is a monad on $\mathcal{C}$, and $H \colon \mathcal{D} \to \mathcal{C}$, we write $H_* \colon [\mathcal{D},\mathcal{D}] \to [\mathcal{D},\mathcal{C}]$ for the functor that postcomposes with $H$. Then  $(H_* \downarrow TH)$ is the category with objects $(P,\alpha)$ where $P \colon \mathcal{D} \to \mathcal{D}$ is an endofunctor and $\alpha \colon HP \to TH$ is a natural transformation, and morphisms $\beta \colon (P,\alpha) \to (P',\alpha')$ are natural transformations $\beta \colon P \to P'$ satisfying $\alpha'H\beta = \alpha$. As suggested above, when $\mathcal{C}=\textbf{Meas}$ and $\mathbb{T}=\mathcal{G}$, we can interpret objects of $(H_* \downarrow \mathcal{G}H)$ as \textit{endofunctors of probability measures}, and its terminal object as the \textit{maximal monad of probability measures}, since it inherits a monad structure from $\mathcal{G}$.
\begin{lemma}
    If $\mathbb{T}$ is a monad on $\mathcal{C}$, $H \colon \mathcal{D} \to \mathcal{C}$ and $(P,\alpha)$ is the terminal object of $(H_* \downarrow TH)$ then $P$ has a unique monad structure making $\alpha$ a Kleisli law.
\end{lemma}
\proof
This immediately follows from 
\begin{enumerate}
    \item Monoids in $(H_* \downarrow TH)$ coincide with monads admitting a Kleisli law into $\mathbb{T}$, since the forgetful functor $U \colon (H_* \downarrow TH) \to [\mathcal{D},\mathcal{D}]$ is strict monoidal
    \item Terminal objects have a unique monoid structure
\end{enumerate} 
However, for clarity, we spell out the proof. First, we have $\mu^\mathbb{T}_H T\alpha \alpha_P \colon HPP \to TH$, so there is a unique map $\mu^\mathbb{P} \colon PP \to P$ satisfying $\alpha H\mu^\mathbb{P} = \mu^\mathbb{T}_H T\alpha \alpha_P$. Also, we have $\eta^\mathbb{T}_H \colon H \to TH$, so there is a unique map $\eta^\mathbb{P} \colon 1_\mathcal{D} \to P$ satisfying $\alpha H \eta^\mathbb{P} = \eta^\mathbb{T}_H$. Note that these are exactly the equations for $\alpha$ to be a Kleisli law (assuming $\mathbb{P}$ is indeed a monad),    and so this is the unique possible structure. Now $\eta^\mathbb{P}_P \colon P \to PP$ satisfies
\begin{align*}
\mu^\mathbb{T}_H T\alpha \alpha_P H\eta^\mathbb{P}_P &= \mu^\mathbb{T}_H T\alpha \eta^\mathbb{T}_{HP} \\
&=\mu^\mathbb{T}_H \eta^\mathbb{T}_{TH}\alpha \\
&=\alpha
\end{align*}
so it is a morphism in $(H_* \downarrow TH)$, and hence $\mu^\mathbb{P}\eta^\mathbb{P}_P = 1_P$, since $(P,\alpha)$ is terminal. Similarly $P\eta^\mathbb{P}$ satisfies
\begin{align*}
\mu^\mathbb{T}_H T\alpha \alpha_P HP\eta^\mathbb{P} &= \mu^\mathbb{T}_H T\alpha TH\eta^\mathbb{P} \alpha \\
&=\mu^\mathbb{T}_H T\eta^\mathbb{T}_H\alpha \\
&=\alpha
\end{align*}
so we also have $\mu^\mathbb{P} P\eta^\mathbb{P} = 1_P$. Finally,  $\mu^\mathbb{T}_HT\mu^\mathbb{T}_H TT\alpha T\alpha_P \alpha_{PP} \colon HPPP \to TH$ defines an object in $(H_* \downarrow TH)$ and $\mu^\mathbb{P}_P$ satisfies
\begin{align*}
\mu^\mathbb{T}_H T\alpha \alpha_P H\mu^\mathbb{P}_P &= \mu^\mathbb{T}_H T\alpha \mu^\mathbb{T}_{HP} T\alpha_P \alpha_{PP} \\
&=\mu^\mathbb{T}_H T\mu^\mathbb{T}_H TT\alpha T\alpha_P \alpha_{PP}
\end{align*}
but $P\mu^\mathbb{P}$ also satisfies
\begin{align*}
\mu^\mathbb{T}_H T\alpha \alpha_P HP\mu^\mathbb{P} &= \mu^\mathbb{T}_H T\alpha TH\mu^\mathbb{P} \alpha_{PP} \\
&=\mu^\mathbb{T}_H T\mu^\mathbb{T}_H TT\alpha T\alpha_P \alpha_{PP}
\end{align*} so both are morphisms in $(H_* \downarrow TH)$, and hence associativity holds.
\endproof
\begin{example}
\textbf{Countable expectation monad on Set.} \label{example4.5}

If $H \colon \mathcal{D}\to \mathcal{C}$ has a right adjoint $U$, then $(H^* \downarrow TH) \cong [\mathcal{D},\mathcal{D}]/UTH$, which has a terminal object given by $UTH$. Hence, let $H$ be the functor $D \colon \textbf{Set} \to \textbf{Meas}$ which assigns a set its discrete measurable space, as in \hyperref[example2.4]{Section 2.4}. Then $\mathcal{E}_c$ and the map $D\mathcal{E}_c \to \mathcal{G}D$ with components the identity function form the terminal object in $(D_* \downarrow \mathcal{G}D)$.
\end{example}
 Note that terminal objects of $(H_* \downarrow TH)$ are exactly the right Kan lifts of $TH$ against $H$, written $(P,\alpha)=\text{Rift}_H TH $.  
 \subsection{Codensity monads and liftings}
 In this section, we will show that the probability monads from \hyperref[section2]{Section 2} are terminal in full subcategories of $(H_* \downarrow \mathcal{G}H)$, for appropriate $H$ in each case. A limited form of this universal property was observed in \cite{van_breugel_metric_nodate}.
 \begin{example}~\cite{van_breugel_metric_nodate} \label{exmp4.7}
 If $H \colon \textbf{KMet} \to \textbf{Meas}$ assigns the Borel $\sigma${-}algebra to a compact metric space, then $H\mathcal{K}=\mathcal{G}H$, and $1_{H\mathcal{K}} \colon H\mathcal{K} \to \mathcal{G}H$ is a Kleisli law (or explicitly, $\eta^\mathcal{G}_H=H\eta^\mathcal{K}$ and $\mu^\mathcal{G}_H=H\mu^\mathcal{K}$). Furthermore, $\text{ev}_{\{0\}} \colon K2 \to [0,1]$ is an isomorphism in \textbf{KMet}. Let $\mathbb{T}$ be a monad on \textbf{KMet} such that $HT=\mathcal{G}H$, $1_{HT}$ is a Kleisli law, and $\text{ev}_{\{0\}} \colon T2 \to [0,1]$ is an isomorphism. Then the map $T \to \mathcal{K}$ whose components are the identity function (but is not necessarily an identity morphism in $\textbf{KMet}$) is a monad morphism.
 \end{example}
 Note that if $T \colon \textbf{KMet} \to \textbf{KMet}$ is any endofunctor satisfying $HT=\mathcal{G}H$ then the underlying set of $TX$ is the set of Borel probability measures on $X$ and for any $f \colon Y \to X$, $Tf$ is given by sending measures to their pushforward along $f$. Also, since $H$ is the identity on morphisms, the condition that $1_{HT}$ should be a Kleisli law uniquely determines a monad structure on $T$, and this monad structure exists iff the functions $\eta_X \colon X \to TX$ and $\mu_X \colon TTX \to TX$ are non{-}expansive for each $X$. Hence, the result in \cite{van_breugel_metric_nodate} shows that  $(\mathcal{K},1_{H\mathcal{K}})$ is terminal in a full subcategory of $(H_* \downarrow \mathcal{G}H)$ whose objects are $(T,1_{HT})$ such that $\eta_X \colon X \to TX$ and $\mu_X \colon TTX \to TX$ are non{-}expansive for each $X$ and such that $\text{ev}_{\{0\}} \colon T2 \to [0,1]$ is an isomorphism in $\textbf{KMet}$.
 
 As remarked in \cite{van_breugel_metric_nodate}, this universal property shows that the Kantorovich distance is maximal among a class of admissible metrics for the topology of weak convergence of probability measures. In \hyperref[example4.12]{Example 4.12}, we use the codensity presentation of $\mathcal{K}$ to show that $\mathcal{K}$ is terminal in a wider subcategory of $(H_* \downarrow TH)$. Hence, we can conclude that the Kantorovich distance is maximal among a larger class of possible metrics for probability measures. 
\begin{example}
\textbf{Radon monad on KHaus.}\label{example4.8}

Let $H \colon \textbf{KHaus} \to \textbf{Meas}$ send a compact Hausdorff space to its Baire measurable structure, the $\sigma$-algebra generated by \emph{zero sets} $f^{-1}(\{0\})$ for $f\in C(X,[0,1])$ (see \hyperref[appenA]{Appendix A}). As shown in \hyperref[corA.4]{Corollary A.4}, there is an isomorphism $\text{res} \colon H\mathcal{R} \to \mathcal{G}H$, with $\text{res}_{[0,1]}$ an identity morphism. However, it is not clear that $(\mathcal{R},\text{res})$ is terminal in $(H_* \downarrow \mathcal{G}H)$. If $\beta \colon HS \to \mathcal{G}H$ factors through $\text{res}$, we must have $\beta_{[0,1]} = H\lambda$ for some (necessarily unique) $\lambda$. 
\end{example}
 This motivates the following result, which we will use in \hyperref[example4.11]{Example 4.11} to show that this condition is the only impediment to terminality.
\begin{theorem}
\label{thm4.9}
Consider the following lax square
\[\begin{tikzcd}
	{\mathcal{F}} & {\mathcal{D}} \\
	{\mathcal{E}} & {\mathcal{C}}
	\arrow["F"', from=1-1, to=2-1]
	\arrow["I", from=1-1, to=1-2]
	\arrow["K", from=1-2, to=2-2]
	\arrow["H"', from=2-1, to=2-2]
	\arrow["\gamma"', shorten <=0.3cm, shorten >=0.3cm, Rightarrow, from=2-1, to=1-2]
\end{tikzcd}\]
Let $(\mathbb{T},\varepsilon)$ be the codensity monad of $K$ and suppose that
\begin{enumerate}[label=(\roman*)]
        \item $H$ and $I$ are faithful, and $\gamma$ is monic 
        \item $(TH,\varepsilon_I T\gamma)= \mathrm{Ran}_{F} KI$
    \end{enumerate}
Then for $P\colon \mathcal{E}\to \mathcal{E}$, the following are equivalent
\begin{enumerate}
    \item There is a $\rho \colon PF\to F$ making $(\mathbb{P},\rho)$ the codensity monad of $F$ 
    \item There is a map $\alpha \colon HP \to TH$ making $(P,\alpha)$ the terminal object in the full monoidal subcategory of $(H_* \downarrow TH)$ of objects $(S,\beta)$ where $\varepsilon_I T\gamma \beta_F=\gamma H\nu$ for a (necessarily unique) $\nu \colon SF \to F$
\end{enumerate} Furthermore, the monad structures induced on $P$ by these universal properties coincide.
\end{theorem}
\proof
We first verify that this is indeed a monoidal subcategory (by which we mean that the inclusion should be strict monoidal). We note that if $(S,\beta)$ and $(R,\delta)$ are objects in $(H_* \downarrow TH)$, their monoidal product is given by $(SR,\mu^\mathbb{T}_H T\delta \beta_R)$, and if we have $\varepsilon_I T\gamma \beta_F=\gamma H\nu$ and $\varepsilon_I T\gamma \delta_F=\gamma H\kappa$ then
\begin{align*}
\varepsilon_I T\gamma \mu^\mathbb{T}_{HF} T\delta_F \beta_{RF}  &=\varepsilon_I \mu^\mathbb{T}_{KI} TT\gamma T\delta_F \beta_{RF}\tag{naturality} \\
&=\varepsilon_I T[\varepsilon_I T\gamma \delta_F] \beta_{RF}\tag{definition of $\mu^\mathbb{T}$} \\
&=\varepsilon_I T\gamma TH\kappa \beta_{RF} \\
&=\varepsilon_I T\gamma \beta_F  HS\kappa \tag{naturality} \\
&=\gamma H(\nu S\kappa) \\
\end{align*}
The unit in $(H_* \downarrow TH)$ is given by $(1_\mathcal{E},\eta^\mathbb{T}_H)$ and $\varepsilon_I T\gamma\eta^\mathbb{T}_{HF}=\gamma$.

Now, given any $S \colon \mathcal{E} \to \mathcal{E}$ and $\nu \colon SF \to F$ there is a unique map $\beta \colon HS \to TH$ such that $\varepsilon_I T\gamma \beta_F=\gamma H\nu$, since $(TH,\varepsilon_I T\gamma)= \mathrm{Ran}_{F} KI$. Also, if $\lambda \colon S \to R$, $\kappa \colon RF \to F$, and $\delta \colon HR\to TH$ is such that $\varepsilon_I T\gamma \delta_F=\gamma H\kappa$ then 
\begin{align*}
\nu=\kappa \lambda_F & \iff \gamma H\nu=\gamma H\kappa H\lambda_F \\
& \iff  \varepsilon_I T\gamma \beta_F = \varepsilon_I T\gamma \delta_F H\lambda_F \\
& \iff \beta =\delta H\lambda
\end{align*}
  Let $F^* \colon [\mathcal{E},\mathcal{E}]\to [\mathcal{F},\mathcal{E}]$ be the precomposition functor. There is a full, faithful, and injective on objects functor $L \colon (F^* \downarrow F) \to (H_* \downarrow TH)$ that sends $(S,\nu)$ to $(S,\beta)$, and acts as the identity on morphisms. Furthermore, since the monoidal product of $(S,\nu)$, and $(R,\kappa)$ in $(F^* \downarrow F)$ is $(SR,\nu S \kappa)$ and the unit is $(1_\mathcal{E},1_F)$, $L$ is strict monoidal. Now, a terminal object $(P,\nu)$ in $(F^* \downarrow F)$  is sent by $L$ to an object $(P,\alpha)$, which is then the terminal object of the image of $L$. Furthermore, since $L$ is full and faithful, if $L(P,\nu)=(P,\alpha)$ is terminal in the image of $L$, then clearly $(P,\nu)$ is terminal in $(F^* \downarrow F)$.  Finally, since $L$ acts as the identity on morphisms and is strict monoidal, the monad structures on $P$ coincide.
\endproof
One immediate consequence of this result is that we can construct Kleisli laws of probability monads from their codensity presentations.
\begin{corollary}
    In the setting of \hyperref[thm4.9]{Theorem 4.9}, if $(\mathbb{P},\rho)$ is the codensity monad of $F$, then there is a Kleisli law $\alpha \colon HP\to TH$.
\end{corollary}
\begin{example}\label{example4.11} \textbf{Radon monad on KHaus continued.}

We apply \hyperref[thm4.9]{Theorem 4.9} to show that, for the Baire functor $H$, the Radon monad is terminal in the full subcategory of $(H_* \downarrow \mathcal{G}H)$ of objects $(S,\beta)$ such that $\beta_{[0,1]} \colon HS[0,1] \to \mathcal{G}[0,1]$ is continuous (so that $\beta_{[0,1]}=H\lambda$ for some $\lambda$). We let $F=K_\mathcal{R}, K=K_\mathcal{G}$, and $I, H$ be the functors in \hyperref[example4.8]{Example 4.8}. Since the isomorphism $\text{res} \colon H\mathcal{R} \to \mathcal{G}H$ from \hyperref[corA.4]{Corollary A.4} is such that $\text{res}_{\mathcal{R}X}$ is the identity for finite discrete spaces $X$, we let $\gamma=1_{HK_\mathcal{R}}$. In \hyperref[appenA]{Appendix A} it is verified that $\mathcal{G}H=\textnormal{Ran}_{K_\mathcal{R}}K_\mathcal{G}I$, and the remaining hypotheses of \hyperref[thm4.9]{Theorem 4.9} also hold. Since $\mathcal{R}$ is the codensity monad of $K_\mathcal{R}$, \hyperref[thm4.9]{Theorem 4.9} shows that $(\mathcal{R}, \text{res})$ is terminal in the full subcategory of $(H_*\downarrow \mathcal{G}H)$ on $(S,\beta)$ such that $\mu^\mathcal{G}_X\beta_{\mathcal{R}X}$ is continuous at finite discrete space $X$ (we are using $X$ to both denote a discrete measurable and topological space). Now, if $X$ is a finite discrete space, $\mu^\mathcal{G}_X \beta_{K_\mathcal{R}X}$ is continuous iff $\mathcal{G}\mathbbm{1}_{\{x\}}\mu^\mathcal{G}_X \beta_{K_\mathcal{R}X}$ is continuous for each $x\in X$, where $\mathbbm{1}_{\{x\}} \colon X \to \{0,1\}$ is the indicator function. Furthermore, by naturality we have $\mathcal{G}\mathbbm{1}_{\{x\}}\mu^\mathcal{G}_X \beta_{K_\mathcal{R}X} = \mu^\mathcal{G}_{\{0,1\}}\beta_{[0,1]}HS\mathcal{R}\mathbbm{1}_{\{x\}}$. Hence, since $\mu^\mathcal{G}_{\{0,1\}}$ is continuous, it is sufficient for $\beta_{[0,1]}$ to be continuous, and we saw this is also necessary in \hyperref[example4.8]{Example 4.8}.
\end{example}
\begin{example} \label{example4.12} \textbf{Kantorovich monad on \textbf{KMet}.}

 We apply \hyperref[thm4.9]{Theorem 4.9} to show that, for the Borel functor $H \colon \textbf{KMet}\to \textbf{Meas}$, the Kantorovich monad is terminal in the full subcategory of $(H_* \downarrow \mathcal{G}H)$ of objects $(S,\beta)$ such that $\beta_{[0,1]} \colon HS[0,1] \to \mathcal{G}[0,1]$ is non-expansive. 
 
 First we note that the functor $H\colon \mathbf{KMet}\to \mathbf{Meas}$ that assigns a compact metric space its Borel $\sigma$-algebra preserves any limits that exist in $\mathbf{KMet}$. We prove this by first noting that the functor $T\colon \mathbf{KMet}\to \mathbf{KHaus}$ which assigns a metric space its induced topology preserves any limits that exist, since in the commutative triangle
 \[\begin{tikzcd}
	{\mathbf{KMet}} & {\mathbf{KHaus}} \\
	& {\mathbf{Set}}
	\arrow["T", from=1-1, to=1-2]
	\arrow["{\mathbf{KMet}(1,-)}"', from=1-1, to=2-2]
	\arrow["{\mathbf{KHaus}(1,-)}", from=1-2, to=2-2]
\end{tikzcd}\]
the functor $\mathbf{KMet}(1,-)$ preserves any limits that exist, and $\mathbf{KHaus}(1,-)$ creates all limits in $\mathbf{KHaus}$. Hence $T$ must preserve any limits that exist, and combining this with~\cite[Theorem 7.1.1]{dudley_real_2002}, which states that the Baire $\sigma$-algebra of a metric space coincides with its Borel $\sigma$-algebra and \hyperref[thmA.3]{Theorem A.3}, it follows that $H$ preserves limits that exist. Thus, if $I \colon \textbf{FinStoch} \to \textbf{cStoch}$ is the inclusion,  $H\mathcal{K}X=\lim_J HK_\mathcal{K}= \lim_J K_\mathcal{G}I$ and thus $\mathcal{G}H=\text{Ran}_{K_\mathcal{K}}K_{\mathcal{G}}I$.

 Hence, we can apply the same argument as in \hyperref[example4.11]{Example 4.11}. This generalises the terminality result in \cite{van_breugel_metric_nodate}, discussed in \hyperref[exmp4.7]{Example 4.7}, which also required a monad structure on objects $(S,1_{HS})$ such that $\mu,\eta$ are non-expansive. Intuitively, this shows that if we can prove that an inequality of metrics of probability measures holds for $X=[0,1]$, we can derive it for any compact metric space $X$.
 
As an illustrative example, we now recall that if $X$ is a metric space and $p,q$ are Borel probability measures on $X$, their Prokhorov distance is defined by
 \[\mathrm{d}_{PX}(p,q)= \inf \{\alpha > 0 : \forall A \in \mathcal{B}X  \;  p(A) \leq q(A_\alpha)+\alpha \}\]
 where $A_\alpha =\{x : \mathrm{d}_X(x,A) \leq \alpha \}$ when $A \neq \emptyset$ and $\emptyset_\alpha =\emptyset$~\cite{parthasarathy_probability_1967,billingsley_convergence_1999}. Let $PX$ be the space of Borel probability measures with the Prokhorov metric, then as shown in \cite{billingsley_convergence_1999}, if $X$ is separable then $PX$ has the topology of weak convergence; it has the initial topology for $\text{ev}_h \colon PX \to [0,1]$ for $h \in C(X,[0,1])$. In particular, if $(X,\text{d}_X)$ is a compact metric space, so is $(PX,\mathrm{d}_{PX})$. Furthermore, if $f \colon X \to Y$ is non{-}expansive, we can define $Pf \colon PX \to PY$ to send $p$ to its pushforward along $f$, defined by $Pf(p)(A)=p(f^{{-}1}(A))$ Then $Pf$ is also non{-}expansive since
 \begin{align*}
  \mathrm{d}_{PX}(p,q) & = \inf \{\alpha > 0 : \forall A \in \mathcal{B}X  \;  p(A) \leq q(A_\alpha)+\alpha \} \\
& \geq \inf \{\alpha > 0 : \forall B \in \mathcal{B}Y  \;  p(f^{{-}1}(B)) \leq q(f^{{-}1}(B)_\alpha)+\alpha \} \\
& \geq \inf \{\alpha > 0 : \forall B \in \mathcal{B}Y \;  p(f^{{-}1}(B)) \leq q(f^{{-}1}(B_\alpha))+\alpha \} \\
& = \mathrm{d}_{PY}(Pf(p),Pf(q))
 \end{align*}
In the penultimate line we used that $f$ is non-expansive to conclude that $f^{{-}1}(B)_\alpha \subseteq f^{{-}1}(B_\alpha)$. Hence, $P$ defines an endofunctor on \textbf{KMet} with $HP=\mathcal{G}H$. In \cite{van_breugel_metric_nodate} it is suggested that we might use this in conjunction with a terminality result for $\mathcal{K}$ to show that $\mathrm{d}_{PX} \leq \mathrm{d}_{\mathcal{K}X}$. However, our result does not apply (and this inequality cannot hold) since the identity function $\text{id}_{[0,1]} \colon P[0,1] \to \mathcal{K}[0,1]$ fails to be non{-}expansive. For example, if we define $p=\frac{3}{5}\eta_{[0,1]}(0)+\frac{2}{5}\eta_{[0,1]}(1)$ and $q=\frac{3}{5}\eta_{[0,1]}(\frac{1}{2})+\frac{2}{5}\eta_{[0,1]}(1)$, then $\mathrm{d}_{\mathcal{K}X}(p,q)=\frac{3}{10}$, but $\mathrm{d}_{PX}(p,q) = \frac{1}{2}$. We do have the inequality  $\mathrm{d}_{PX}^2 \leq \mathrm{d}_{\mathcal{K}X}$ for any complete separable metric space $X$~\cite[Corollary 4.3]{huber_robust_1981}. 
\end{example}
\begin{example} \textbf{Expectation monad on Set.}

Suppose $H \colon \textbf{Set} \to \textbf{Meas}$ assigns the discrete measurable structure, $I=1_{\textbf{cStoch}}$ and $\gamma=\text{id} \colon HK_{\mathcal{E}_c} \to K_\mathcal{G}$ has components the identity function. Then the subcategory specified in \hyperref[thm4.9]{Theorem 4.9} is the entire category $(H_* \downarrow TH)$.
\end{example} 
\subsection{Codensity presentations from the lifting property}
 We conclude this section by showing another consequence of \hyperref[thm4.9]{Theorem 4.9}; that we can derive a Kleisli codensity presentation for a universal lifting of a codensity monad. In particular, this gives a new method to derive codensity presentations of probability monads that have a universal property as a lifting of the Giry monad. Suppose now that $(P,\alpha)$ is terminal in $(H_* \downarrow TH)$. By \hyperref[prop4.3]{Proposition 4.3}, if $\alpha$ is monic and $H$ is faithful, there is a faithful functor $\Bar{H} \colon \mathcal{D}_\mathbb{P} \to \mathcal{C}_\mathbb{T}$. We will show that, under suitable conditions, $\mathbb{P}$ and $\mathbb{T}$ share a Kleisli codensity presentation over a small subcategory of $\mathcal{D}_\mathbb{P}$, via the functor $\Bar{H}$.
 \begin{corollary} \label{cor4.15}
   Suppose $\mathbb{T}$ is the codensity monad of $K \colon \mathcal{D} \to \mathcal{C}$, and $(P,\alpha)$ is the terminal object in $(H_* \downarrow TH)$ for a faithful functor $H \colon \mathcal{E} \to \mathcal{C}$. If there is a faithful functor $J \colon \mathcal{D} \to \mathcal{E}_\mathbb{P}$ such that $\alpha_J$ is pointwise monic and $K=G_\mathbb{T}\Bar{H}I$, where $\Bar{H}$ is the lifting induced by $\alpha$, then $\mathbb{P}$ and $\mathbb{T}$ both have a Kleisli codensity presentation over $\mathcal{D}$. 
\end{corollary}
First, $\bar{H}J$ is faithful by \hyperref[prop4.3]{Proposition 4.3} and so $\mathbb{T}$ has a Kleisli codensity presentation over $\mathcal{D}$. Now, in \hyperref[thm4.9]{Theorem 4.9}, we let $F=G_\mathbb{P}J$, $I=1_\mathcal{D}$ and $\gamma=\alpha_J$ (at least component{-}wise). Then it suffices to show that $(P,\alpha)$ is an object of the full subcategory specified in \hyperref[thm4.9]{Theorem 4.9}, since it is terminal in all of $(H_* \downarrow TH)$. But this is demonstrated by the Kleisli law equation $\alpha H\mu^\mathbb{P}=\mu^\mathbb{T}_H T\alpha \alpha_P$. Hence, by \hyperref[thm4.9]{Theorem 4.9}, $\mathbb{P}$ has a Kleisli codensity presentation over $\mathcal{D}$.
\endproof
\begin{example}
\textbf{Countable expectation monad on Set continued.} \label{exmp4.16}

Let $H$ be the functor $D \colon \textbf{Set} \to \textbf{Meas}$ which assigns a set its discrete $\sigma$-algebra, as in \hyperref[example4.5]{Example 4.5}. One can verify that all the conditions of \hyperref[cor4.15]{Corollary 4.15} hold, and so $\mathcal{E}_c$ has a Kleisli codensity presentation over $\textbf{cStoch}$. In this case, $D$ has a right adjoint $U$, and so this can also be directly shown via the isomorphisms \[U(\text{Ran}_{K_{\mathcal{G}}} K_\mathcal{G})D \cong \text{Ran}_U (\text{Ran}_{K_\mathcal{G}} UK_\mathcal{G}) \cong \text{Ran}_{UK_\mathcal{G}} UK_\mathcal{G}\]
and using that $UK_\mathcal{G}=K_{\mathcal{E}_c}$. But \hyperref[cor4.15]{Corollary 4.15} can be applied more generally, when the functor $H$ does not necessarily have a right adjoint. 
\end{example}

\section{Commutativity of codensity monads and probability bimeasures} \label{section5}
 The definition of a codensity monad in \hyperref[section3]{Section 3} can be generalised to any 2{-}category, by requiring that the diagram in \hyperref[definition3.2]{Definition 3.2} is a Kan extension in that 2{-}category. In this section, we begin by giving sufficient conditions for a codensity monad to receive a lax monoidal structure, making it the codensity monad in the 2{-}category \textbf{MonCat} of monoidal categories, lax monoidal functors, and monoidal natural transformations. In \hyperref[prop5.6]{Proposition 5.6}, we give a simplified condition for the probability monads considered in \hyperref[section2]{Section 2}. Verifying this condition will rely on the theory of probability bimeasures, and we provide a basic overview of this in \hyperref[subsection5.8]{Section 5.8}. We then show that the Radon monad and the restriction of the Giry monad to standard Borel spaces satisfy the condition in \hyperref[prop5.6]{Proposition 5.6}, but the other monads from \hyperref[section2]{Section 2} do not. In fact, the Radon monad and restriction of the Giry monad satisfy much stronger conditions, making them \textit{exactly pointwise monoidal codensity monads}, and we study these conditions in \hyperref[subsection5.18]{Section 5.18}. The main result of this subsection is \hyperref[thm5.23]{Theorem 5.23}, which provides a characterisation of the monoidal codensity monads that are exactly pointwise, and gives a description of their Kleisli category.
 
 \subsection{Conditions for commutativity}   
 
  We begin by providing sufficient conditions for a codensity monad to be lax monoidal. These are similar to conditions found in \cite{fritz_criterion_nodate} for lifting left Kan extensions to \textbf{MonCat}, and are closely related to the general results in \cite{koudenburg_algebraic_2015,weber_algebraic_2016} on lifting Kan extensions to 2{-}categories of algebras for 2{-}monads.

\begin{theorem}

\label{thm5.2}
Suppose $(K,\kappa,\iota) \colon \mathcal{D} \to \mathcal{C}$ is a lax monoidal functor and that $(\mathbb{T},\varepsilon)$ is the codensity monad of $K$. Then if 
\begin{enumerate}[label=(\roman*)]
\item $\varepsilon_I T\iota \colon TI \to KI $ is monic
\item  $\varepsilon_{\text{{-}}\otimes \text{{-}}}T\kappa \colon T(K{-}\otimes K{-}) \to K({-}\otimes {-})$ makes $T({-}\otimes {-})$ the right Kan extension of $K({-}\otimes {-})$ along $K\times K$
\item $\varepsilon_{(\text{{-}}\otimes \text{{-}})\otimes \text{{-}}}T(\kappa_{\text{{-}}\otimes \text{{-}},\text{{-}}} \, \kappa \otimes 1_K) \colon T((K{-}\otimes K{-})\otimes K{-}) \to K(({-}\otimes {-})\otimes {-})$ makes $T(({-}\otimes {-})\otimes {-})$ the right Kan extension of $K(({-}\otimes {-})\otimes {-})$ along $(K\times K)\times K$
\end{enumerate}
then there is a unique monoidal structure $\chi$ on $T$ making $(\mathbb{T},\varepsilon)$ the codensity monad of $K$ in $\mathbf{MonCat}$. \\
If furthermore $K$ is symmetric monoidal, then so is $T$, and $(\mathbb{T},\varepsilon)$ is the codensity monad of $K$ in $\mathbf{SMonCat}$.
\end{theorem}
\proof
First, by $(ii)$, $\kappa(\varepsilon\otimes \varepsilon) \colon TK{-}\otimes TK{-} \to  K({-}\otimes {-})$ has a unique factorisation $\chi \colon T{-}\otimes T{-} \to T({-}\otimes {-})$ such that $\varepsilon_{\text{{-}}\otimes \text{{-}}}T\kappa \, \chi_{K\times K}= \kappa(\varepsilon\otimes \varepsilon)$. We claim that $(T,\chi, \eta_I)$ is a monoidal structure on $T$. Since $K$ is lax monoidal, $T\lambda \, \chi_{I,\text{{-}}}(\eta_I\otimes 1_T) \colon I\otimes T{-} \to T$ satisfies
\begin{align*}
    \varepsilon \,  T\lambda_K \,  \chi_{I,K\text{{-}}}(\eta_I\otimes 1_{TK}) & = \varepsilon \,  TK\lambda \, T\kappa_{I,\text{{-}}}T(\iota\otimes 1_K)\chi_{I,K\text{{-}}}(\eta_I\otimes 1_{TK}) \\
    & =K\lambda \, \varepsilon_{I,\text{{-}}}T\kappa_{I,\text{{-}}} \, \chi_{KI,K\text{{-}}}(T\iota \, \eta_I\otimes 1_{TK}) \tag{naturality}\\
    & =K\lambda \, \kappa_{I,\text{{-}}}((\varepsilon_I T\iota \, \eta_I)\otimes \varepsilon) \tag{definition of $\chi$}\\
    & =K\lambda \, \kappa_{I,\text{{-}}}((\varepsilon_I \eta_{KI}  \iota)\otimes  \varepsilon) \tag{naturality}\\
    & = K\lambda \, \kappa_{I,\text{{-}}}(\iota \otimes 1_K)(1_I \otimes  \varepsilon) \tag{definition of $\eta$}\\
    & = \lambda_K(1_I \otimes  \varepsilon) \tag{$K$ lax monoidal}\\
    & =\varepsilon \, \lambda  _{TK} \tag{naturality}
\end{align*}
and so, by the uniqueness of factorisation, $T\lambda  \chi_{I,\text{{-}}}(\eta_I\otimes 1_T)=\lambda  _T$. Dually, the other unitor square commutes. This calculation also shows $\varepsilon_I T\iota \, \eta_I=\iota$, which is one of the conditions for $\varepsilon$ to be monoidal. 

Next, $T\alpha \, \chi_{\text{{-}},\text{{-}}\otimes \text{{-}}}(1_T\otimes \chi) \colon T{-}\otimes (T{-}\otimes T{-}) \to T(({-}\otimes {-})\otimes {-})$ satisfies
\begin{align*}
&\varepsilon_{(\text{{-}}\otimes \text{{-}})\otimes \text{{-}}}T[\kappa_{\text{{-}}\otimes \text{{-}},\text{{-}}} \, (\kappa\otimes 1_K)]T\alpha_{K,K,K}\chi_{K\text{{-}},K\text{{-}}\otimes K\text{{-}}}(1_{TK}\otimes \chi_{K,K}) \\
 &= \varepsilon_{(\text{{-}}\otimes \text{{-}})\otimes \text{{-}}}T[K\alpha  \kappa_{\text{{-}},\text{{-}}\otimes \text{{-}}}(1_K\otimes \kappa)]\chi_{K\text{{-}},K\text{{-}}\otimes K\text{{-}}}(1_{TK}\otimes \chi_{K,K}) \tag{$K$ 
 monoidal}\\
 &= K\alpha \, \varepsilon_{\text{{-}}\otimes \text{{-}}\otimes (\text{{-}})}T\kappa_{\text{{-}},\text{{-}}\otimes \text{{-}}} \chi_{K\text{{-}},K(\text{{-}}\otimes \text{{-}})}(1_{TK}\otimes T\kappa \, \chi_{K,K}) \tag{naturality} \\
 &=K\alpha \, \kappa_{\text{{-}},\text{{-}}\otimes \text{{-}}} (\varepsilon \otimes  \varepsilon_{\text{{-}}\otimes \text{{-}}}) \, T\kappa \, \chi_{K,K}) \tag{definition of $\chi$} \\
 &=K\alpha \, \kappa_{\text{{-}},\text{{-}}\otimes \text{{-}}} (1_K \otimes  \kappa) (\varepsilon \otimes  (\varepsilon \otimes  \varepsilon)) \tag{definition of $\chi$} \\
&=\kappa_{\text{{-}}\otimes \text{{-}},\text{{-}}} (\kappa \otimes  1_K) \alpha_{K,K,K} \varepsilon \otimes  (\varepsilon \otimes  \varepsilon)) \tag{$K$ monoidal} \\
 &=\kappa_{\text{{-}}\otimes \text{{-}},\text{{-}}} (\kappa 
 (\varepsilon \otimes  \varepsilon)\otimes  \varepsilon)\alpha_{TK,TK,TK} \tag{naturality} \\
 &=\kappa_{\text{{-}}\otimes \text{{-}},\text{{-}}}(\varepsilon_{\text{{-}}\otimes \text{{-}}}\otimes  \varepsilon)(T\kappa \, \chi_{K,K} \otimes  1_{TK})\alpha_{TK,TK,TK} \tag{definition of $\chi$} \\
 &=\varepsilon_{(\text{{-}}\otimes \text{{-}})\otimes \text{{-}}}T\kappa_{\text{{-}}\otimes \text{{-}},\text{{-}}} \chi_{K(\text{{-}}\otimes \text{{-}}),K\text{{-}}}(T\kappa \, \chi_{K,K} \otimes  1_{TK})\alpha_{TK,TK,TK} \tag{definition of $\chi$} \\
 &=\varepsilon_{(\text{{-}}\otimes \text{{-}})\otimes \text{{-}}}T[\kappa_{\text{{-}}\otimes \text{{-}},\text{{-}}}(\kappa \otimes  1_{K})] \chi_{K\text{{-}}\otimes K\text{{-}},K\text{{-}}}( \, \chi_{K,K} \otimes  1_{TK})\alpha_{TK,TK,TK} \tag{naturality}
\end{align*}
Thus, by the uniqueness of factorisation, $T$ is lax monoidal.
Also, $\varepsilon \colon TK \to K$ is monoidal since the coherence relating to $\kappa$ and $\chi$ is the defining equation for $\chi$. Finally, suppose $\delta \colon GK \to K $ is a monoidal natural transformation, where $(G,\xi, \gamma)$ is a lax monoidal endofunctor on $\mathcal{C}$. Then there is a unique natural transformation $\varphi \colon G \to T$ satisfying $\varepsilon \, \varphi_K = \delta$. We show that $\varphi$ monoidal.
\begin{align*}
\varepsilon_{A \otimes B}T\kappa_{A,B}\varphi_{KA \otimes KB}\xi_{KA,KB} &= \varepsilon_{A \otimes B}\varphi_{K(A \otimes B)}G\kappa_{A,B}\xi_{KA,KB} \tag{naturality}\\
& = \delta_{A \otimes B}G\kappa_{A,B}\xi_{KA,KB} \tag{definition of $\varphi$} \\
& = \kappa_{A,B}(\delta_A \otimes  \delta_B) \tag{$\delta$  monoidal} \\
& = \kappa_{A,B}(\varepsilon_A\varphi_{KA} \otimes  \varepsilon_B \varphi_{KB}) \tag{definition of $\varphi$} \\
& = \varepsilon_{A \otimes B}T\kappa_{A,B}\chi_{KA,KB}(\varphi_{KA} \otimes \varphi_{KB}) \tag{definition of $\chi$} \\
\end{align*}
Hence, $\varphi_{A \otimes B}\xi_{A,B}=\chi_{A,B}(\varphi_A \otimes \varphi_B)$ follows by the uniqueness of factorisation of the Kan extension in $(ii)$.  
Next,
\begin{align*}
\varepsilon_I T\iota \, \eta_I &= \iota \\
& = \delta_I G\iota \, \gamma \tag{$\delta$ monoidal} \\
& = \varepsilon_I \varphi_{KI} G\iota \, \gamma \tag{definition of $\alpha $} \\
& = \varepsilon_I T\iota \, \varphi_I \gamma \tag{naturality} \\
\end{align*}
Hence, since $\varepsilon_I T\iota \ $ is monic by assumption, it follows that $\varphi$ is monoidal.

If $\chi'$ is another lax monoidal structure making $(\mathbb{T},\varepsilon)$ the codensity monad in \textbf{MonCat}, then since $\varepsilon$ is monoidal we must have $\varepsilon_{\text{{-}}\otimes \text{{-}}}T\kappa\chi'_{K\times K}= \kappa(\varepsilon\otimes \varepsilon)$. Hence by the uniqueness of factorisation of the Kan extension in $(ii)$, $\chi=\chi'$.

Finally, if $K$ is also symmetric monoidal, then $T\gamma_{B,A} \chi_{B,A} \colon TB \otimes  TA \to T(A \otimes  B)$ satisfies 
\begin{align*}
\varepsilon_{A\otimes B} T \kappa_{A,B} T\gamma_{KB,KA} \chi_{KB,KA} &= \varepsilon_{A\otimes B} TK\gamma_{B,A} T\kappa_{B,A} \chi_{KB,KA} \\
&=F\gamma_{B,A} \varepsilon_{B\otimes B} T\kappa_{B,A} \chi_{KB,KA} \tag{naturality} \\
&=K\gamma_{B,A} \kappa_{B,A}(\varepsilon_B \otimes  \varepsilon_A)  \tag{by definition of $\chi$}\\
&= \kappa_{A,B} \gamma_{KB,KA} (\varepsilon_B \otimes  \varepsilon_A) \tag{$K$ symmetry} \\
&= \kappa_{A,B} (\varepsilon_A \otimes  \varepsilon_B) \gamma_{TKB,TKA} \tag{naturality}\\\
&=\varepsilon_{A\otimes B} T \kappa_{A,B} \chi_{KA,KB} \gamma_{TKB,TKA} \tag{definition of $\chi$}
\end{align*}
Hence, $T\gamma_{B,A}\chi_{B,A} = \chi_{A,B} \gamma_{TB,TA}$ by the uniqueness of factorisation through the Kan extension in $(ii)$ and so $T$ is symmetric monoidal.
\endproof
If $\mathbb{T}$ and the Kan extension in condition $(ii)$ are pointwise, then for $h \colon A \to KX$, $k \colon B \to KY$ we will write $\text{ev}_{h,k} \colon T(A\otimes  B) \to K(X \otimes  Y)$ for the cone maps of the limits in $(ii)$. We then have 
\begin{align*}
\text{ev}_{1_{KX},1_{KY}} &= \varepsilon_{X \otimes  Y} T\kappa_{X,Y}\\
&=\text{ev}_{\kappa_{X,Y}}
\end{align*}
Then, for $f \colon A' \to A$ and $ g \colon B' \to B $, we have  $\text{ev}_{h,k}T(f\otimes  g)=\text{ev}_{hf,kg}$ and so
\begin{align*}
\text{ev}_{h,k} &= \text{ev}_{\kappa_{X,Y}} T(h\otimes  k) \\
&=\text{ev}_{\kappa_{X,Y} (h\otimes  k)}
\end{align*}
The lax monoidal strength $\chi$ is the unique map satisfying $\text{ev}_{\kappa_{X,Y} (h\otimes  k)} \chi_{A,B} = \kappa_{X,Y}(\text{ev}_h\otimes  \text{ev}_k)$. Note that in this setting, condition $(i)$ states that the map $\text{ev}_\iota$ should be monic.
\begin{remark}  If condition $(i)$ holds and $\iota$ is an isomorphism, then since $\varepsilon_I T\iota \eta_I =\iota$, $\varepsilon_I T\iota$ is an isomorphism, and hence $\eta_I$ is too. Next, note that if the conditions of \hyperref[prop3.11]{Proposition 3.11} hold, then $\varepsilon_I T\iota$ is an isomorphism, so is certainly monic, and we can drop condition $(i)$. Finally, if $\mathcal{D}$ is a monoidal category which is a subcategory of $\mathcal{C}_\mathbb{T}$, $\mathbb{T}$ has a Kleisli codensity presentation over $\mathcal{D}$, and $\iota \colon I \to G_\mathbb{T}I=TI$ is given by $\eta_I$, then that $\varepsilon_I T\iota=1_{TI}$ is one of the unit identities for a monad, and so in this case we may also drop condition $(i)$.
\end{remark}
 It will be convenient to have simpler assumptions than those of \hyperref[thm5.2]{Theorem 5.2} for the settings in which we apply it. We have already seen that in many situations we can drop condition $(i)$. The following results will show that for the probability monads from \hyperref[section2]{Section 2}, condition $(iii)$ can also be dropped.
 \begin{definition}[Concrete pointed (monoidal) category]
     A locally small category $\mathcal{C}$ with a terminal object $1$ is concrete pointed if $\mathcal{C}(1,{-})$ is faithful. A (symmetric) monoidal category $\mathcal{C}$ is concrete pointed if furthermore $\mathcal{C}(1,{-})$ is strong (symmetric) monoidal with respect to the Cartesian structure on \textnormal{\textbf{Set}}.
 \end{definition}
 The categories in \hyperref[section2]{Section 2} are all concrete pointed as symmetric monoidal categories. It is shown in \cite{moggi_notions_1991} that any endofunctor on a concrete pointed monoidal category has a unique strength, if it exists. The result provides a mild generalisation, by dropping the associativity coherence axiom for strengths.
\begin{lemma}[\cite{moggi_notions_1991}]
\label{lem5.5}
If $\mathcal{C}$ is a concrete pointed monoidal category, and $F:\mathcal{C} \to \mathcal{C}$ has maps $\theta_{A,B} \colon A \otimes F B \to F(A \otimes B)$  natural in $A,B$ such that $F\lambda_A \theta_{1,A} = \lambda_{FA}$, then $\theta$ is a unique left strength for $F$.
\end{lemma}
\proof
 First, if $\mathcal{C}$ is concrete pointed monoidal, then the unit $1$ is terminal and any map $f \colon 1 \to A\otimes B$ factors as $(f_1 \otimes f_2) \lambda_1^{{-}1} \colon 1 \to 1\otimes 1 \to A \otimes B$ for unique $f_1, f_2$.  Now if $a \colon 1 \to A $, $h \colon 1 \to FB$, then 
\begin{align*}
\mathcal{C}(1,\theta_{A,B})((a\otimes h) \lambda_1^{{-}1}) &= \mathcal{C}(1,\theta_{A,B}(a\otimes 1_{FB}))((1_1\otimes h) \lambda_1^{{-}1}) \\
&= \mathcal{C}(1,F(a\otimes 1_B) \theta_{1,B})((1_1\otimes h) \lambda_1^{{-}1}) \tag{by naturality}\\
&= \mathcal{C}(1,F[(a \otimes 1_B)\lambda_B^{{-}1}]\lambda_{FB})((1_1\otimes h)\lambda_1^{{-}1}) \tag{by hypothesis} \\
&= \mathcal{C}(1,F[(a \otimes 1_B)\lambda_B^{{-}1}])(h)
\end{align*} and hence since $\mathcal{C}(1,{-})$ is faithful, $\theta$ is uniquely defined. One can then verify that this does indeed satisfy the associativity coherence law, using the fact that $\mathcal{C}(1,{-})$ is strong monoidal. Hence, $\theta$ is the unique left strength as desired.
\endproof
 This lemma shows that for concrete pointed monoidal categories, the associativity coherence law for strengths is determined by the other axioms. We thus obtain associativity of (symmetric) monoidal monads for free in these categories.
\begin{lemma}
\label{lem5.6}
If $\mathcal{C}$ is a concrete pointed (symmetric) monoidal category and $\mathbb{T}$ is a monad on $\mathcal{C}$ with $\chi_{A,B} \colon TA \otimes TB \to T(A \otimes B)$ natural in $A,B$ such that the following diagrams commute
\begin{mathpar}
\begin{tikzcd}
	{I\otimes TA} & TA \\
	{TI\otimes TA} & {T(I\otimes A)}
	\arrow["{\lambda_{TA}}", from=1-1, to=1-2]
	\arrow["{ \iota\otimes1_{TA}}"', from=1-1, to=2-1]
	\arrow["{\chi_{I,A}}"', from=2-1, to=2-2]
	\arrow["{T\lambda_A}"', from=2-2, to=1-2]
\end{tikzcd} \and
\begin{tikzcd}
	{TA\otimes I} & TA \\
	{TA\otimes TI} & {T(A\otimes I)}
	\arrow["{\rho_{TA}}", from=1-1, to=1-2]
	\arrow["{1_{TA}\otimes \iota}"', from=1-1, to=2-1]
	\arrow["{\chi_{A,I}}"', from=2-1, to=2-2]
	\arrow["{T\rho_A}"', from=2-2, to=1-2]
\end{tikzcd}
\and
\begin{tikzcd}
	{A\otimes B} & {TA\otimes TB} \\
	& {T(A\otimes B)}
	\arrow["{\eta_A\otimes \eta_B}", from=1-1, to=1-2]
	\arrow["{\eta_{A\otimes B}}"', from=1-1, to=2-2]
	\arrow["{\chi_{A,B}}", from=1-2, to=2-2]
\end{tikzcd}
\and
\begin{tikzcd}
	{TTA\otimes TTB} & {T(TA\otimes TB)} & {TT(A\otimes B)} \\
	{TA\otimes TB} && {T(A\otimes B)}
	\arrow["{\chi_{TA,TB}}", from=1-1, to=1-2]
	\arrow["{\mu_A \otimes \mu _B}"', from=1-1, to=2-1]
	\arrow["{T\chi_{A,B}}", from=1-2, to=1-3]
	\arrow["{\mu_{A\otimes B}}"', from=1-3, to=2-3]
	\arrow["{\chi_{A,B}}"', from=2-1, to=2-3]
\end{tikzcd}
\end{mathpar}
Then $\mathbb{T}$ is a (symmetric) lax monoidal monad. In particular, the coherence laws for $\chi$ with the associator and symmetry additionally hold.
\end{lemma}
\proof
Define $\theta_{A,B}=\chi_{A,B}(\eta_A \otimes 1_{TB}) \colon A \otimes TB \to T(A \otimes B)$ then $\theta$ satisfies the conditions of \hyperref[lem5.5]{Lemma 5.5}  and so $\theta$ is the unique left strength of $T$. Similarly $\vartheta_{A,B}=\chi_{A,B}(1_{TA} \otimes \eta_B)$ is the unique right strength, and $\mathbb{T}$ is a (symmetrically) bistrong monad. Hence, 
\begin{align*}
\chi_{A,B} &= \chi_{A,B}(\mu_A T\eta_{A}\otimes \mu_B T\eta_B) \\
&=\mu_{A\otimes B}T\chi_{A,B}\chi_{TA,TB}(T\eta_{A}\otimes \eta_{TB}) \tag{by coherence for $\mu$}\\
&=\mu_{A\otimes B}T\chi_{A,B}T(\eta_A \otimes 1_TB)\chi_{A,TB}(1_{TA} \otimes \eta_{TB}) \\
&=\mu_{A\otimes B}T\theta_{A,B}\vartheta_{A,TB}
\end{align*}
But this is the derived monoidal strength from the left and right strengths, which will satisfy the associativity coherence law, and so $\chi$ satisfies it too. The symmetry coherence law also follows when $\mathcal{C}$ is concrete pointed symmetric monoidal.
\endproof
 In particular, in \hyperref[thm5.2]{Theorem 5.2}, if $\mathcal{C}$ is a concrete pointed (symmetric) monoidal category, then we can drop condition $(iii)$ since this was only used to prove the associativity of $\chi$. All the conditions of \hyperref[lem5.6]{Lemma 5.6} hold since the existence of $\chi$ and its unit laws were derived from condition $(ii)$. That the Kan extension factorisation map $\varphi \colon GK \to K$ is monoidal was also derived from conditions $(i)$ and $(ii)$, and did not rely on the associativity of the lax monoidal strength $\xi$ of $G$. Hence, the proof of \hyperref[thm5.2]{Theorem 5.2} shows that under conditions $(i)$ and $(ii)$, $(\mathbb{T},\chi)$ is the codensity monad of $(K,\kappa,\iota)$ in the 2{-}category of monoidal categories, lax monoidal functors that may not satisfy the associativity coherence law, and monoidal natural transformations between such functors. Hence, $\mu$ and $\eta$ are monoidal, and by \hyperref[lem5.6]{Lemma 5.6} $\mathbb{T}$ is (symmetric) monoidal,even if $K$ was not assumed to be symmetric.
\begin{proposition}
\label{prop5.6}
Suppose $\mathcal{C}$ is a concrete pointed (symmetric) monoidal category and $\mathbb{T}$ is a monad on $\mathcal{C}$ with a Kleisli codensity presentation over $\mathcal{D}$. Let $K_\mathcal{D}$ be the restriction of the Kleisli forgetful functor $G_\mathbb{T} \colon \mathcal{C}_\mathbb{T} \to \mathcal{C}$ to $\mathcal{D}$. If there is a monoidal structure on $\mathcal{D}$, and a lax monoidal structure $(\kappa, \eta_I)$ on $K_\mathcal{D}$, such that condition $(ii)$ of \hyperref[thm5.2]{Theorem 5.2} is satisfied, then $\mathbb{T}$ is a (symmetric) lax monoidal monad.
\end{proposition}
\subsection{Spaces of probability $k${-}polymeasures} \label{subsection5.8}  

 We now give a review of the basic theory of bimeasures, and more generally $k${-}polymeasures. We use this to show that for the probability monads studied in \hyperref[section2]{Section 2}, the limit condition in \hyperref[prop5.6]{Proposition 5.6} requires that the space of bimeasures coincides with the space of bivariate probability measures. We show that this holds for the Radon monad, but fails for the Giry monad unless it is restricted to standard Borel spaces.
\begin{definition}[Probability $k${-}polymeasure]
    A \textit{probability bimeasure} $\gamma$ on measurable spaces $X,Y$ is a function $\gamma \colon \mathcal{B}X \times \mathcal{B}Y \to [0,1]$ that is separately countably additive in each component and satisfies $\gamma (X,Y)=1$. More generally, a \textit{probability $k${-}polymeasure} on measurable spaces $X_1,\dots,X_k$ is a function $\gamma \colon \prod_{i=1}^k \mathcal{B}X_i \to [0,1]$ that is countably additive in each component and satisfies $\gamma(X_1,\dots,X_k)=1$. 
\end{definition}
\begin{example} \label{example5.10}
    If $p$ is a probability measure on $\prod_{i=1}^k X_i$, then $p$ defines a probability $k${-}polymeasure $\gamma$ on $X_1,\dots,X_k$ by $\gamma(A_1,\dots,A_k)=p( \prod_{i=1}^k A_i)$. In this case, we say the polymeasure $\gamma$ \textit{extends to} $p$.
\end{example}
\begin{example} \label{example5.11}
    If $p$ is a probability measure on $X$ and $q$ is a probability measure on $Y$, then $\gamma(A,B)=p(A)\cdot q(B)$ defines a bimeasure on $X,Y$. More generally, if $\gamma_1,\dots,\gamma_n$ are polymeasures of arity $k_1,\dots,k_n$, there is a $(k_1+\dots+k_n)${-}polymeasure given by $ \prod_{i=1}^n\gamma_i(A_{i1},\dots,A_{ik_i})$
\end{example}
\begin{example}~\cite{karni_extension_1990}
\label{example5.12}

This example shows that there is a probability bimeasure on measurable spaces $X,Y$ which does not extend to a measure on $X\times Y$. Let $m$ be the Lebesgue measure on $[0,1]$ and pick two disjoint sets $X,Y \subseteq [0,1]$ of first category with outer Lebesgue measure 1. Now we let  $\mathcal{B}X= \{A \cap X : A \in \mathcal{B}[0,1] \}$ and define $\mathcal{B}Y$ analogously. Define $\gamma$ on $X,Y$ by $\gamma (A',B')=m(A\cap B)$ when $A'=A \cap X$, $B'=B \cap Y$, $A,B \in \mathcal{B}[0,1]$. If we let $E_i^n=[\frac{i}{2^n},\frac{i+1}{2^n}]\cap X$ and $F_i^n$ is defined analogously, then if  $Q_n=\bigcup_{i=1}^{2^n{-}1}E_i^n \times F_i^n $, we have $\sum_{i=1}^{2^n{-}1} \gamma(E_i^n ,F_i^n)=1$ and $m(Q_n)=1/2^n$. Also $\bigcap_{n \in \mathbb{N}}Q_n =\emptyset$ so $\gamma$ cannot extend to a $\sigma${-}additive measure.
\end{example}
 While not all probability bimeasures extend to a measure, Theorem 2.8 of \cite{karni_extension_1990} shows that under weak conditions on measurable spaces $X$ and $Y$, probability bimeasures on $X,Y$ will extend to probability measures on $X\times Y$. In particular, well-behaved $k${-}polymeasures on compact Hausdorff spaces will always extend to measures on the product.
\begin{definition}[Radon probability $k${-}polymeasure]
If $\gamma$ is a Borel probability $k${-}poly\-measure on compact Hausdorff spaces $X_1,\dots,X_k$ (a probability $k${-}polymeasure on the Borel $\sigma${-}algebras), we say that $\gamma$ is Radon if it is Radon separately in each variable. This means that for each $i$, $\gamma(A_1,\dots,A_i,\dots,A_k)=\sup \{\gamma(A_1,\dots,K_i,\dots,A_k) : K_i\subseteq A_i \text{ is compact} \}$.
\end{definition}

\begin{proposition} \label{prop5.14} Every Radon probablity $k${-}polymeasure on compact Hausdorff spaces $X_1,\dots,X_k$  extends to a Radon probability measure on $\prod_{i=1}^k X_i$.
\end{proposition}
\proof
The bimeasure case is \cite[Corollary 2.9]{karni_extension_1990}; see also the historical article \cite{marczewski_remarks_1953}. The $k${-}polymeasure case follows from \cite[Theorem 3.3]{bombal_integral_2001}.
\endproof
If $\gamma$ is a $k${-}polymeasure on $X_1,\dots,X_k$ and $f_i \colon X_i \to [0,1]$ are measurable maps for $i=1,\dots,k$, there is a standard notion of integration with respect to $\gamma$, denoted $\int (f_1,\dots,f_k) \; \text{d}\gamma$. This is defined analogously to integration with respect to a probability measure. It is defined such that it is multilinear, preserves monotone limits in each variable and $\int (\mathbbm{1}_{A_1},\dots,\mathbbm{1}_{A_k}) \; \text{d}\gamma =\gamma(A_1,\dots,A_k)$~\cite{dobrakov_multilinear_1999}.

We can define measurable spaces of $k${-}polymeasures in analogy to the space of probability measures defined in \hyperref[example2.3]{Section 2.3}. If $X_1,\dots,X_k$ are measurable spaces, we define $\mathcal{G}_k(X_1,\dots,X_k)$ to be the set of $k${-}polymeasures on $X_1,\dots,X_k$ with the coarsest $\sigma${-}algebra such that the function $\text{ev}_{A_1,\dots,A_k} \colon \mathcal{G}_k(X_1,\dots,X_k) \to [0,1]$ defined by $\gamma \mapsto \gamma(A_1,\dots,A_k)$ is measurable for each $(A_1,\dots,A_k) \in \prod_{i=1}^k \mathcal{B}X_i $. Furthermore, for $f_i \colon X_i \to [0,1]$, $i=1,\dots,k$, we can define $\text{ev}_{f_1,\dots,f_k} \colon \mathcal{G}_k(X_1,\dots,X_k) \to [0,1]$ by $\text{ev}_{f_1,\dots,f_k}(\gamma)=\int (f_1,\dots,f_k) \; \text{d}\gamma$. Then the $\sigma${-}algebra on $\mathcal{G}_k(X_1,\dots,X_k)$ is equivalently the coarsest such that $\text{ev}_{f_1,\dots,f_k}$ is measurable for each measurable $f_1,\dots,f_k$. Note that $\mathcal{G}_1=\mathcal{G}$  and $\mathcal{G}_0$ is the terminal object in \textbf{Meas}.

We could also define compact Hausdorff spaces of Radon $k${-}polymeasures, but the following example shows that this is not necessary, in addition to verifying that the Radon monad satisfies the conditions of \hyperref[prop5.6]{Proposition 5.6}.
\begin{example}\label{example5.15}
\textbf{Radon monad on KHaus.}

Let $K=K_\mathcal{R}$ and $\mathcal{D}=\textbf{FinStoch}$ as in \hyperref[example3.7]{Example 3.7}. $\textbf{KHaus}$ is concrete pointed monoidal, and since multiplication is continuous, it makes $K_\mathcal{R}$ lax monoidal. Hence, to verify the conditions of \hyperref[prop5.6]{Proposition 5.6}, it remains to show that for any compact Hausdorff spaces $A$ and $B$ we have $\mathcal{R}(A\times B)=\lim_{(A\downarrow K_\mathcal{R}) \times (B \downarrow K_\mathcal{R})} K_{\mathcal{R}}(U_A\times U_B$).

Let $\gamma_X \colon X \times \{0,1\} \to X+1$ be the function defined by $\gamma_X(x,0)=x$, $\gamma_X(x,1)=\bot$ for any countable set $X$. If $(C,\tau_{h,k})$ is a cone over $(A\downarrow K_\mathcal{R}) \times (B \downarrow K_\mathcal{R})$, then for a fixed $k\colon B \to K_\mathcal{R}\{0,1\}$, we have that $(C,\mathcal{R}\gamma_X\tau_{h,k})$ forms a cone over $K_\mathcal{R}(U_A+1) \colon (A\downarrow K_\mathcal{R}) \to \textbf{KHaus}$. Hence, by \hyperref[prop3.12]{Proposition 3.12}, there is a unique map $p_k \colon C \to \mathcal{R}(A+1)$ such that $\text{ev}_{\{(0,0)\}} \tau_{h,k}(c)= \int \text{ev}_{\{0\}}h \, \text{d}p_k(c)$ for any $h \colon B \to K_\mathcal{R}\{0,1\}$. This extends to a continuous positive linear functional $C(A) \to \mathbb{R}$. By the Banach{-}Steinhaus theorem, the map $B_c \colon C(A)\times C(B) \to \mathbb{R}$ given by the linear extension of $B_c(h,k)=\text{ev}_{\{(0,0)\}}\tau_{\text{ev}^{{-}1}_{\{0\}}h,\text{ev}^{{-}1}_{\{0\}}k}(c)$ for $h \colon A \to [0,1]$, $k \colon B \to [0,1]$ defines a continuous bilinear functional. By \hyperref[prop5.14]{Proposition 5.14}, there is a unique Radon measure $p(c)$ on $A \times B$ such that $B_c(h,k)= \int h \cdot k \; \text{d}p(c)$.

We now show the function $p \colon C \to \mathcal{R}(A \times B)$ defined by this unique measure is continuous. First, $\text{ev}_{h\cdot k} p=\text{ev}_{\{(0,0)\}} \tau_{h,k}$ is continuous. Now, define $\mathcal{R}_2(A,B)$ to be the space of Radon probability bimeasures on $A,B$ with the coarsest topology such that the function $\text{ev}_{h,k} \colon \mathcal{R}_2(A,B) \to [0,1]$ (which sends $\gamma \mapsto \int (h,k) \; \text{d}\gamma$) is continuous for each continuous $h \colon A \to [0,1]$ and $k \colon B \to [0,1]$. Then \hyperref[prop5.14]{Proposition 5.14} shows there is a continuous bijection from $\mathcal{R}(A \times B)$ to the set of bimeasures $\mathcal{R}_2(A,B)$. Furthermore, the space $\mathcal{R}_2(A,B)$ can be seen to be a closed subspace of the closed unit ball of $(C(A)\otimes_\pi C(B))^*$ with the weak $*${-}topology, where $\otimes_\pi$ denotes the projective tensor product \cite[Theorem 2]{bowers_representation_2015}. Hence, by the Banach{-}Alaoglu theorem, $\mathcal{R}_2(A,B)$ is a compact Hausdorff space, and so the continuous bijection defined by \hyperref[prop5.14]{Proposition 5.14} is a homeomorphism. Thus, the topology on $\mathcal{R}(A \times B)$ is the coarsest such that $\text{ev}_{h\cdot k}$ is continuous for each $h,k$ as above, and so $p$ is continuous.

Finally we have, for $h \colon A \to K_\mathcal{R} X$, $k \colon B \to K_\mathcal{R} Y$, that 
\begin{align*}
\text{ev}_{h,k}p(c)(x,y) & = \int \text{ev}_{\{x\}}h \cdot \text{ev}_{\{y\}}k \; \text{d}p(c) \\
 & =\int \text{ev}_{\{0\}}\mathcal{R} \mathbbm{1}_{\{x \}} h \cdot \text{ev}_{\{0\}}\mathcal{R} \mathbbm{1}_{\{y \}} k \; \text{d}p(c) \\
 & = \text{ev}_{(0,0)} \tau_{\mathcal{R} \mathbbm{1}_{\{x \}} h,  \mathcal{R} \mathbbm{1}_{\{y \} k}}(c) \\
 & = \tau_{h,k}(c)(x,y)
\end{align*}
Additionally, if $\tau_{h,k}=\text{ev}_{h,k} p'$ then we must have $\text{ev}_{(0,0)}\tau_{h,k}(c)= \int \text{ev}_{\{0\}} h \cdot \text{ev}_{\{0\}} k \; \text{d}p'(c)$ when $X=Y=\{0,1\}$. Thus, by the uniqueness of the Radon measure above, $p'=p$. Hence, $\mathcal{R}(A\times B)$ is the desired limit, and so $\mathcal{R}$ is lax monoidal.

An identical argument shows that $\mathcal{R}(\prod_{i=1}^k A_i) =\lim_{(A_i \downarrow K)}K_\mathcal{R}(\prod_{i=1}^k U_{A_i})$. \endproof
\end{example}
\begin{example} \label{example5.16}
\textbf{Giry monad on \textbf{Meas} and \textbf{BorelMeas}.}\\
Consider the diagram $K_\mathcal{G}(U_A \times U_B) \colon (A \downarrow K_\mathcal{G}) \times (B \downarrow K_\mathcal{G}) \to \textbf{Meas}$ from \hyperref[thm5.2]{Theorem 5.2}. When $h \colon A \to K_\mathcal{G}X$, $k \colon B \to K_\mathcal{G}Y$, the maps $\text{ev}_{h,k}$ (defined by $\text{ev}_{h,k}(\gamma)(x,y)=\int (\text{ev}_{\{x \}} h,\text{ev}_{\{y \}} k) \; \text{d}\gamma $) form a cone over $(A \downarrow K_\mathcal{G}) \times (B \downarrow K_\mathcal{G})$ with apex $\mathcal{G}_2(A,B)$. A similar argument to the one in \hyperref[example5.15]{Example 5.15} shows that this is the limit cone.  \hyperref[example5.12]{Example 5.12} demonstrates that the spaces $\mathcal{G}_2(A,B)$ and $\mathcal{G}(A \times B)$ are not generally isomorphic, and so condition $(ii)$ of \hyperref[thm5.2]{Theorem 5.2} fails in this instance. This argument also shows that $\mathcal{G}_k(A_1,\dots,A_k)=\lim_{\prod_{i=1}^k(A_i \downarrow K_\mathcal{G})}K_\mathcal{G}(\prod_{i=1}^kU_{A_i})$.

Let $\textbf{BorelMeas}$ be the full subcategory of $\textbf{Meas}$ whose objects are standard Borel spaces~\cite{bogachev_measure_2007}. These are measurable spaces $X$ such that there is a admissible metric $\mathrm{d}_X$ on $X$ (so that $\mathcal{B}X$ is the Borel $\sigma${-}algebra of $(X,\mathrm{d}_X)$) making $X$ a complete separable metric space. If $X$ is any standard Borel space, and $\mathrm{d}_X$ is an admissible metric on $X$, then the Prokhorov metric is an admissible metric for $\mathcal{G}X$ and makes it a complete and separable metric space~\cite[Theorem 6.8]{billingsley_convergence_1999}. Hence, there is a monad $\mathcal{G}_B$ on \textbf{BorelMeas} given by the restriction of $\mathcal{G}$~\cite{giry_categorical_1982}. All countable discrete measurable spaces are standard Borel, so we can define $K_{\mathcal{G}_B} \colon \textbf{cStoch} \to \textbf{BorelMeas}$ in the same way as $K_\mathcal{G}$, and then $\mathcal{G}_B$ is the pointwise codensity monad of $K_{\mathcal{G}_B}$. If $A_1,\dots,A_k$ are standard Borel spaces, then the morphism of cones $\mathcal{G}(\prod_{i=1}^k A_i) \to \mathcal{G}_k(A_1,\dots,A_k)$ is an isomorphism by applying Kuratowski's theorem, which gives a classification of standard Borel spaces~\cite[Theorem 13.1.1]{dudley_real_2002}.
\end{example}
 These examples show that, for a probability monad $P$, the limit $\lim_{\prod_{i=1}^k (A_i\downarrow K_P)} K_P$ $(\bigotimes_{i=1}^k U_{A_i})$ often coincides with the space of probability $k${-}polymeasures. When this space coincides with the space of probability measures on $\bigotimes_{i=1}^k A_i$, we can conclude from \hyperref[thm5.2]{Theorem 5.2} that $P$ is a monoidal codensity monad. These are not necessary conditions, since we can directly verify that the Giry monad on \textbf{Meas} is the codensity monad of $K_\mathcal{G}$ in \textbf{SMonCat}, despite not satisfying this condition. We will characterise when these conditions hold for monoidal codensity monads in  \hyperref[subsection5.18]{Section 5.18}. 
 In our examples so far, whether these spaces coincide has been closely related to the extension of probability $k${-}polymeasures to probability measures, but this is not the only condition that needs to be met. In the case of the Kantorovich monad, every probability bimeasure extends to a probability measure, but the metric we obtain on the space of probability bimeasures is strictly smaller than the metric on the space of joint probability measures, and so these spaces are not isomorphic in \textbf{KMet}. In the case of the expectation monad, bimeasures can lift to measures, but not necessarily uniquely. Suppose $\kappa$ is a measurable cardinal, and $\mathcal{U}$ is a free $\kappa${-}complete ultrafilter on $\kappa$. Then, as \hyperref[example2.4]{Section 2.4} showed, the ultrafilter $\mathcal{U}\times \mathcal{U}$ defined on the boolean algebra generated by the rectangles in $\kappa \times \kappa$ extends to $\mathcal{P}(\kappa \times \kappa)$, but not uniquely. 
\begin{remark}
 The limit formula $\mathcal{G}_k(A_1,\dots,A_k)=\lim_{\prod_{i=1}^k(A_i \downarrow K_\mathcal{G})}K_\mathcal{G}(\prod_{i=1}^kU_{A_i})$ exhibits the functors $\mathcal{G}_k \colon \textbf{Meas}^k \to \textbf{Meas}$ as a pointwise right Kan extension, where $\mathcal{G}_k(f_1,\dots,f_k)$ is the unique map satisfying $\text{ev}_{h_1,\dots,h_k} \mathcal{G}_k(f_1,\dots,f_k)=\text{ev}_{h_1f_1,\dots,h_kf_k}$. Explicitly, we have $\mathcal{G}_k(f_1,\dots,f_k)(\gamma)(A_1,\dots,A_k)=\gamma(f_1^{{-}1}(A_1),\dots,f_k^{{-}1}(A_k))$. Together these functors organise into a functor $\mathcal{G} \colon \textbf{Meas}^* \to \textbf{Meas}$. Here $\textbf{Meas}^*$ is the free strict monoidal category on \textbf{Meas}, with objects finite lists $[X_1,\dots,X_k]$ of measurable space $X_i$  and morphisms finite lists $[f_1,\dots,f_k] \colon [X_1,\dots,X_k] \to [Y_1,\dots,Y_k]$ of measurable maps $f_i \colon X_i \to Y_i$. As noted in \hyperref[example5.10]{Example 5.10} and \hyperref[example5.11]{Example 5.11}, the functor $\mathcal{G}$ has additional structure. For example \hyperref[example5.11]{Example 5.11} gives maps:
 \[\prod_{i=1}^n\mathcal{G}[X_{i1},\dots,X_{i k_i}] \to \mathcal{G}[X_{11},\dots,X_{n k_n}]\] In fact, the additional structure of $\mathcal{G}$ can be completely captured by unit maps
\[\eta_{[X_1,\dots,X_n]} \colon \prod_{i=1}^n X_i \to \mathcal{G}[X_1,\dots,X_n]\]
satisfying $\text{ev}_{h_1,\dots,h_n}\eta_{[X_1,\dots,X_n]}(x_1,\dots,x_n)=\prod_{i=1}^n h_i(x_i)$, and multiplication maps
\[\mu_{[[X_{11},\dots,X_{1 k_1}],\dots,[X_{n1},\dots,X_{n k_n}]]} \colon \mathcal{G}[\mathcal{G}[X_{11},\dots,X_{1 k_1}],\dots,\mathcal{G}[X_{n1},\dots,X_{n k_n}]) \to \mathcal{G}[X_{11},\dots,X_{n k_n}]\]
satisfying $\text{ev}_{h_{11},\dots,h_{n k_n}}\mu_{[[X_{11},\dots,X_{1 k_1}],\dots,[X_{n1},\dots,X_{n k_n}]]}=\text{ev}_{\text{ev}_{h_{11},\dots,h_{1 k_1}},\dots,\text{ev}_{h_{n 1},\dots,h_{n k_n}}}$. This gives the functor $\mathcal{G}$ a structure we call a $*${-}monad, which we intend to elaborate further on in future work. We will say here that any lax monoidal functor from a small category to a complete category induces a $*${-}monad in a similar manner.
\end{remark}
\subsection{Exactly pointwise monoidal codensity monads} \label{subsection5.18}   

 We saw in \hyperref[example5.15]{Example 5.15} that the Radon monad satisfies
 \[\mathcal{R}(\prod_{i=1}^k X_i) = \lim_{\prod_{i=1}^k (X_i \downarrow K_\mathcal{R})} K_\mathcal{R}(\prod_{i=1}^k U_{X_i})\] and in \hyperref[example5.16]{Example 5.16} that the Giry monad on standard Borel spaces satisfies a similar condition. An informal view is that a pointwise codensity presentation of a monad $\mathbb{T}$ is a \textit{convenient} limit presentation of the objects $TX$. We can therefore interpret this condition to say that $\mathcal{R}(\prod_{i=1}^k X_i)$ has a convenient limit expression, and hence that $\mathcal{R}$ is a pointwise codensity monad in $\textbf{MonCat}$ in a strong sense.
\begin{definition}[Exactly pointwise monoidal codensity monad] 
    If $(K,\kappa,\iota)$ is a lax monoidal functor and $(\mathbb{T},\varepsilon)$ is its pointwise codensity monad, then $\mathbb{T}$ is \textit{exactly pointwise monoidal in degree} $k\geq 1$ if for every $A_1,\dots,A_k$ we have \[(T(\bigotimes^k_{i=1}A_i),\textnormal{ev}_{\zeta_{X_1,\dots,X_k}\bigotimes_{i=1}^k h_i})=\lim_{\prod_{i=1}^k (A_i \downarrow K)}K(\bigotimes_{i=1}^k U_{A_i})\] where $\zeta_{X_1,\dots,X_k} \colon \bigotimes_{i=1}^k KX_i \to K(\bigotimes_{i=1}^k X_i)$ is the map induced by the lax monoidal structure (see  \cite{epstein_functors_1966,malkiewich_coherence_2022}). It is \textit{exactly pointwise monoidal in degree $0$} if $\text{ev}_\iota \colon TI \to KI$ is an isomorphism. Finally, it is \textit{exactly pointwise monoidal} if it is exactly pointwise monoidal in degree $k$ for every $k \geq 0$. 
\end{definition}

 If $(\mathbb{T},\varepsilon)$ is the pointwise codensity monad of a lax monoidal functor $K$, then it is exactly pointwise monoidal in degree $2$ (resp. $3$) if and only if condition $(ii)$ (resp. $(iii)$) of \hyperref[thm5.2]{Theorem 5.2} holds and is a pointwise Kan extension. Furthermore, if it is exactly pointwise in degree $0$, then condition $(i)$ holds. Hence, exactly pointwise monoidal codensity monads are monoidal codensity monads.
 
The condition for $(\mathbb{T},\varepsilon)$ to be exactly pointwise monoidal is closely related to that of an \textit{exact lax morphism} as used in \cite{street_fibrations_1974,koudenburg_algebraic_2015,weber_algebraic_2016}. Suppose 
\[\begin{tikzcd}
	{\mathcal{C}} & {\mathcal{D}} \\
	{\mathcal{E}}
	\arrow["F", from=1-1, to=1-2]
	\arrow[""{name=0, anchor=center, inner sep=0}, "G"', from=1-1, to=2-1]
	\arrow[""{name=1, anchor=center, inner sep=0}, "R", from=1-2, to=2-1]
	\arrow["\varepsilon"', shift right, shorten <=3pt, shorten >=5pt, Rightarrow, from=1, to=0]
\end{tikzcd}\] is a (pointwise) right Kan extension and 
\[\begin{tikzcd}
	{\mathcal{H}} & {\mathcal{F}} \\
	{\mathcal{C}} & {\mathcal{D}}
	\arrow["K", from=1-1, to=1-2]
	\arrow["H"', from=1-1, to=2-1]
	\arrow["\alpha"', shorten <=6pt, shorten >=6pt, Rightarrow, from=1-2, to=2-1]
	\arrow["L", from=1-2, to=2-2]
	\arrow["F", from=2-1, to=2-2]
\end{tikzcd}\]
is a lax square. Then the (pointwise) right Kan extension is \emph{preserved by pasting the square of $\alpha$} if the composite natural transformation is a (pointwise) right Kan extension, which means  $(RL,\varepsilon_H R\alpha)=\mathrm{Ran}_K(GH)$ (and is pointwise). Now, if $K$ is a lax monoidal functor, its pointwise codensity monad $(\mathbb{T},\varepsilon)$ is exactly pointwise monoidal iff the pointwise right Kan extension $(T,\varepsilon)=\text{Ran}_K K$ is preserved by pasting the lax algebra morphism square for $K$.
\[\begin{tikzcd}
	{\mathcal{D}^*} & {\mathcal{C}^*} \\
	{\mathcal{D}} & {\mathcal{C}}
	\arrow["{K^*}", from=1-1, to=1-2]
	\arrow["{m_\mathcal{D}}"', from=1-1, to=2-1]
	\arrow["\zeta"', shorten <=6pt, shorten >=6pt, Rightarrow, from=1-2, to=2-1]
	\arrow["{m_\mathcal{C}}", from=1-2, to=2-2]
	\arrow["K", from=2-1, to=2-2]
\end{tikzcd}\]
 Here, $\mathcal{C}^*$ is the free strict monoidal category on $\mathcal{C}$, whose objects are finite lists $[X_1,\dots,X_k]$ such that $X_i \in \text{ob} \, \mathcal{C}$, and morphisms are finite lists $[f_1,\dots,f_k] \colon [X_1,\dots,X_k] \to [Y_1,\dots,Y_k]$ where $f_i \colon X_i \to Y_i$ is a morphism in $\mathcal{C}$. The functor $m_\mathcal{C}$ is defined on objects by $[X_1,\dots,X_k] \mapsto \bigotimes_{i=1}^k X_i$ and similarly on morphisms, and $\zeta$ is the map induced by the lax monoidal structure of $K$. More generally, if $(K,\kappa,\iota)$ and $(G,\chi,\nu)$ are lax monoidal functors, then we say that a pointwise right Kan extension $\text{Ran}_K G$ in \textbf{Cat} is exactly pointwise monoidal if and only if it is preserved by pasting with the lax algebra morphism square of $K$. In \cite{weber_algebraic_2016} and \cite{koudenburg_algebraic_2015}, $(K,\zeta)$ is said to be \textit{exact} if every pointwise right Kan extension along $K$ is preserved by pasting the lax algebra morphism square of $K$. In our terminology, $(K,\zeta)$ is exact iff every pointwise right Kan extension along $K$ is exactly pointwise monoidal. Our notion of exactly pointwise monoidal codensity monads coincides with the definition of $f$-exactness in \cite{koudenburg_algebraic_2015,koudenburg_algebraic_2013}.

The remainder of the section will build up to \hyperref[thm5.23]{Theorem 5.23}, which gives a characterisation of when a monoidal codensity monad is exactly pointwise. \hyperref[prop3.3]{Proposition 3.3} shows that when $\mathcal{D}$ is small, and $\mathcal{C}$ is locally small and complete, the codensity monad of a functor $K\colon \mathcal{D}\to \mathcal{C}$ has a resolution via the adjunction $K_\circ \dashv  \mathrm{Ran}_{Y_\mathcal{D}}K$. The Kleisli category of the codensity monad $(\mathbb{T},\varepsilon)$ is then isomorphic to the full subcategory of $[\mathcal{D},\textbf{Set}]^{\text{op}}$ on objects $\mathcal{C}(A,K{-})$ for $A \in \text{ob} \; \mathcal{C}$. We will show in \hyperref[thm5.23]{Theorem 5.23} that this picture lifts to \textbf{MonCat} precisely when the codensity monad $(\mathbb{T},\varepsilon)$ is exactly pointwise monoidal. Recall first that if $\mathcal{D}$ is a (symmetric) monoidal category and $F,G \colon \mathcal{D} \to \textbf{Set}$, their Day convolution~\cite{day_closed_1970} is defined by the coend \[F\otimes_{\text{Day}}G\; ({-})=\int^{(X_1,X_2) \in \text{ob} \, \mathcal{D}\times\mathcal{D}} \mathcal{D}(X_1\otimes X_2,{-})\times FX_1 \times GX_2\] Explicitly, \[F\otimes_{\text{Day}}G(X)=\coprod_{X_1,X_2 \in \text{ob} \, \mathcal{D}} \{(x_1,x_2,m) : x_1 \in FX_1, x_2 \in GX_2, m \colon X_1 \otimes X_2 \to X \}/\sim\] where $\sim$ is the smallest equivalence relation such that for every $f \colon X'_1 \to X_1$ and $g \colon X'_2 \to X_2$ in $\mathcal{D}$, $(x_1,x_2,m(f\otimes g))\sim(Ff(x_1),Fg(x_2),m)$. We denote an equivalence class under this relation by $[x_1,x_2,m]$. This monoidal product makes $[\mathcal{D},\textbf{Set}]$ a (symmetric) monoidal category, with unit $\mathcal{D}(I,{-})$. Furthermore, there is an isomorphism $ [\mathcal{D},\textbf{Set}](F\otimes_{\text{Day}}G,H) \cong [\mathcal{D}\times \mathcal{D},\textbf{Set}](F{-}_1 \times G{-}_2,H({-}_1 \otimes {-}_2))$ written $\alpha \mapsto \Bar{\alpha}$ where $\Bar{\alpha}_X([x_1,x_2,m])=Hm \alpha_{X_1,X_2}(x,y)$~\cite[Section 6.2]{loregian_coend_2021}. For the remainder of the section, if $\mathcal{D}$ is a monoidal category and we refer to a monoidal structure on  $[\mathcal{D},\textbf{Set}]$ (or $[\mathcal{D},\textbf{Set}]^\text{op}$) we will be referring to its Day monoidal structure.\\
We say that an adjunction $F \dashv G$ lifts to \textbf{MonCat} if there are lax monoidal structures on $F,G$ such that the unit and counit of the adjunction are monoidal natural transformations. The following result shows that $F$ must be strong monoidal.
\begin{proposition}[\cite{im_universal_1986,kelly_doctrinal_1974}] \label{prop5.20}
Let $(F,\kappa,\iota) \colon \mathcal{C} \to \mathcal{D}$ be a lax monoidal functor. Then $F$ is a left adjoint in $\mathbf{MonCat}$ iff $F\dashv G$ in $\mathbf{Cat}$ and $F$ is strong monoidal. In this case, there is a unique lax monoidal structure on $G$ so that $\eta$ and $\varepsilon$ are monoidal.
\end{proposition}
\proof
    If $F$ has a right adjoint $(G,\chi,\nu)$ in \textbf{MonCat}, then $F\dashv G$ in \textbf{Cat}, and the maps $\varepsilon_{FA\otimes FB}F\chi_{FA,FB} F(\eta_A \otimes \eta_B) \colon F(A\otimes B) \to FA\otimes FB$ and $\varepsilon_I F\nu \colon FI \to I$ provide inverses to $\kappa_{A,B}$ and $\iota$ respectively.

    Conversely, if $F \dashv G$ in \textbf{Cat} and $F$ is strong monoidal, then $\nu=G\iota^{{-}1}\eta_I$ and $\chi_{A,B}=G(\varepsilon_A \otimes \varepsilon_B)G\kappa^{{-}1}_{GA,GB} \eta_{GA\otimes GB}$ define a monoidal structure on $G$ such that $\eta$ and $\varepsilon$ are monoidal, and this is the unique possible structure making $\eta$ monoidal. 
\endproof
Hence, the adjunction $K_\circ \dashv \mathrm{Ran}_{Y_\mathcal{D}}K$ does not immediately lift to \textbf{MonCat}, even if $(\mathbb{T},\varepsilon)$ is exactly pointwise monoidal, since $K_\circ$ is not generally strong monoidal. However, it is always oplax monoidal.
\begin{lemma} \label{lem5.21}
    Suppose $\mathcal{D}$ is small, $\mathcal{C}$ is locally small and complete, and $K \colon \mathcal{D} \to \mathcal{C}$ is a lax monoidal functor, then $K_\circ \colon \mathcal{C} \to [\mathcal{D},\mathbf{Set}]^{\mathrm{op}}$ is oplax monoidal. Furthermore, if  $\xi^{(k)} \colon K_\circ(\bigotimes_{i=1}^k{-}) \to \bigotimes_{i=1}^kK_\circ({-})$ is the map induced by the oplax monoidal structure and $(\mathbb{T},\varepsilon)$ is the codensity monad of $K$, then it is exactly pointwise monoidal in degree $k$ iff $\mathrm{Ran}_{Y_\mathcal{D}}K(\xi^{(k)})$ is an isomorphism.
\end{lemma}
\proof
    Let $\xi^{(0)} \colon \mathcal{D}(I,{-}) \to \mathcal{C}(I,K{-})$ have components defined by $\xi^{(0)}_X(x)=Kx\iota$ and $\xi^{(2)}_{A,B} \colon \mathcal{C}(A,K{-})\otimes \mathcal{C}(B,K{-}) \to \mathcal{C}(A\otimes B, K{-})$ be the natural transformation corresponding to
    $\Bar{\xi}^{(2)}_{A,B}$ defined by $(\Bar{\xi}^{(2)}_{A,B})_{X,Y}(h,k)=\kappa_{X,Y}(h\otimes k)$. Then, it is routine to verify that this defines a lax monoidal structure for $K_\circ^{\text{op}} \colon \mathcal{C}^\text{op} \to [\mathcal{D},\textbf{Set}]$, or equivalently, an oplax monoidal structure for $K_\circ$.
    
    Now, for $k\geq 1$, the functor $H \colon \prod_{i=1}^k (A_i \downarrow K) \to (* \downarrow \bigotimes_{i=1}^k \mathcal{C}(A_i,K{-}))$ which is defined on objects by \[H((X_1,h_1),\dots,(X_k,h_k))=(\bigotimes_{i=1}^k X_i,[\dots[h_1,h_2,1_{X_1\otimes X_2}],\dots,h_k,1_{\bigotimes_{i=1}^k X_i}])\] is cofinal. That is, $(H\downarrow X)$ is (non-empty and) connected for every $X \in \text{ob} (* \downarrow \bigotimes_{i=1}^k \mathcal{C}(A_i,K{-}))$. This can be seen for $k=2$ since a map $m' \colon H((X'_1,x'_1),(X'_2,x'_2)) \to (X,[x_1,x_2,m])$ is exactly an equivalence $(x'_1,x'_2,m')\sim(x_1,x_2,m)$, by the equivalence relation defining  $\mathcal{C}(A,K{-})\otimes \mathcal{C}(A,K{-})(X)$. An argument by induction shows this for larger $k$. For instance, when $k=3$, $[[h_1,h_2,n],h_3,m]=[[h_1,h_2,1_{X_1 \otimes X_2}],h_3, m(n\otimes 1_{X_3})]$.
    
    Hence, we see $\mathrm{Ran}_{Y_\mathcal{D}}K(\xi^{(k)})_{(A_1,\dots,A_k)} \colon T(\bigotimes_{i=1}^kA_i) \to \lim_{\prod_{i=1}^k(A_i \downarrow K)}K(\bigotimes_{i=1}^k U_{A_i})$ is the unique morphism of cones from $(T(\bigotimes^k_{i=1}A_i),\text{ev}_{\zeta_{X_1,\dots,X_k}\bigotimes_{i=1}^k h_i})$ and so is an isomorphism iff this is a limit cone. Finally, $\mathrm{Ran}_{Y_\mathcal{D}}K(\xi^{(0)})= \text{ev}_\iota$ so the result also holds for $k=0$.
\endproof
Now, let $\mathcal{A}$ be the full subcategory of $ [\mathcal{D},\mathbf{Set}]^{\mathrm{op}}$ on objects $\mathcal{C}(A,K{-})$ for $A \in \text{ob} \; \mathcal{C}$. This lemma suggests that if $
\Bar{\mathcal{A}}$ is the smallest full monoidal subcategory of $ [\mathcal{D},\mathbf{Set}]^{\mathrm{op}}$ containing $\mathcal{A}$, then when $(\mathbb{T},\varepsilon)$ is exactly pointwise monoidal, $\mathcal{A}$ is a coreflective subcategory of $\Bar{\mathcal{A}}$. We will see that this indeed holds, and is part of the characterisation of exactly pointwise monoidal codensity monads in \hyperref[thm5.23]{Theorem 5.23}. The following proposition, which is a variant of Day's reflection theorem~\cite{day_note_1973}, will be used to define a monoidal structure on $\mathcal{A}$.
\begin{proposition}[\cite{day_note_1973}] \label{prop5.22}
    If $\mathcal{D}$ is a reflective subcategory of a monoidal category $\mathcal{C}$ and \allowbreak $(L,\eta) \colon \mathcal{C} 
 \to \mathcal{C}$ is the idempotent monad induced by the reflection, then $L$ has a lax monoidal structure iff $L(\eta_A\otimes \eta_B)$ is an isomorphism for all objects $A,B$. Such a monoidal structure on $L$ is unique.
\end{proposition}
\proof
Uniqueness follows from $L(\eta_A\otimes\eta_B)\kappa_{A,B}=\eta_{LA\otimes LB}$. We will now show $\kappa_{A,B}=L(\eta_A\otimes\eta_B)^{{-}1}\eta_{LA\otimes LB}$ defines a lax monoidal structure on $L$ when $L(\eta_A\otimes \eta_B)$ is an isomorphism. First, $L(\eta_A\otimes \eta_B)^{{-}1}\eta_{LA\otimes LB} (\eta_A\otimes \eta_B)=\eta_{A\otimes B}$ follows from naturality, so $\eta$ is monoidal with respect to this structure. Also, we can show that under these conditions $L((\eta_A\otimes \eta_B)\otimes \eta_C)$ is also an isomorphism for any object $C$ since $L((\eta_A\otimes \eta_B)\otimes \eta_C) \kappa_{A\otimes B,C}=\kappa_{LA\otimes LB,LC}(L(\eta_A \otimes \eta_B)\otimes L\eta_C)$. Also, $L\kappa$ is an isomorphism, so $LL((\eta_A\otimes \eta_B)\otimes \eta_C)$ is too, and hence $L((\eta_A\otimes \eta_B)\otimes \eta_C)$ is. Thus, we have:
\begin{align*}
& L\alpha_{A,B,C}L(\eta_A \otimes \eta_{B\otimes C})^{{-}1}\eta_{LA\otimes L(B\otimes C)}(1_{LA}\otimes(L(\eta_B\otimes \eta_C)^{{-}1}\eta_{LB\otimes LC})) \\
 & = L\alpha_{A,B,C}L(\eta_A \otimes \eta_{B\otimes C})^{{-}1}L(1_{LA}\otimes L(\eta_B\otimes \eta_C))^{{-}1}L(1_{LA}\otimes\eta_{LB\otimes LC}))\eta_{LA\otimes (LB \otimes LC)} \\
 & = L\alpha_{A,B,C}L(\eta_A \otimes L(\eta_B\otimes \eta_C)\eta_{B\otimes C})^{{-}1}L(1_{LA}\otimes\eta_{LB\otimes LC}))\eta_{LA\otimes (LB \otimes LC)} \\
 & = L\alpha_{A,B,C}L(\eta_A \otimes (\eta_B \otimes \eta_C))^{{-}1}\eta_{LA \otimes (LB \otimes LC)} \\
 & = L((\eta_A \otimes \eta_B) \otimes \eta_C)^{{-}1}L\alpha_{LA,LB,LC}\eta_{LA \otimes (LB \otimes LC)} \\
 & = L(L(\eta_A \otimes \eta_B)\eta_{A\otimes B}\otimes \eta_C)^{{-}1}L(\eta_{LA\otimes LB}\otimes 1_{LC})\eta_{(LA\otimes LB)\otimes LC}\alpha_{LA,LB,LC} \\
 & = L(\eta_{A\otimes B}\otimes \eta_C)^{{-}1}\eta_{L(A\otimes B)\otimes LC}(L(\eta_A \otimes \eta_B)^{{-}1}\eta_{LA\otimes LB}\otimes 1_{LC})\alpha_{LA,LB,LC}
\end{align*}
Next, we have
\begin{align*}
L\lambda_A L(\eta_I \otimes \eta_A)^{{-}1}\eta_{LI\otimes LA}(\eta_I\otimes 1_{LA}) & = L\lambda_AL(\eta_I \otimes \eta_A)^{{-}1} L((\eta_I\otimes 1_{LA})) \eta_{I\otimes LA} \\ 
& = L\eta_A^{{-}1} L\lambda_{LA} \eta_{I\otimes LA} \\ 
& = \lambda_{LA}
\end{align*}
and dually for $\rho$.

Conversely, if $L$ is monoidal, then $L\kappa_{A,B}L(\eta_A\otimes \eta_B)=L\eta_{A\otimes B}$ so $L(\eta_A\otimes \eta_B)$ is left invertible. Furthermore, we have $L(\eta_A\otimes \eta_B)\kappa_{A,B}=\eta_{LA\otimes LB}$, using $L\eta=\eta_{L}$. Hence, we see that $LL(\eta_A\otimes \eta_B)$ is right invertible and so $L(\eta_A\otimes \eta_B)$ is an isomorphism.
\endproof
Under the conditions of \hyperref[prop5.22]{Proposition 5.22}, $\mathcal{D}$ has a monoidal structure such that the inclusion $i: \mathcal{D} \to \mathcal{C}$ is a right adjoint in \textbf{MonCat}. Explicitly, this is defined by $A\otimes_\mathcal{D}B=L(A\otimes_\mathcal{C} B)$ and $I_\mathcal{D}=LI_\mathcal{C}$, and $\eta$ gives the lax monoidal structure of $i$. When $\mathcal{D}$ has this monoidal structure, we say $\mathcal{D}$ is a lax monoidal reflective subcategory of $\mathcal{C}$. Dually, $\mathcal{D}^{\text{op}}$ is an oplax monoidal coreflective subcategory of $\mathcal{E}^{\text{op}}$. Now, we can give the characterisation of exactly pointwise monoidal codensity monads. First, we note that if $\mathcal{B}$ is any full subcategory of $[\mathcal{D},\mathbf{Set}]^{\mathrm{op}}$ containing $\mathcal{A}$, there is an adjunction $K_\circ |^\mathcal{B} \dashv \text{Ran}_{Y_\mathcal{D}}K |_\mathcal{B}$ where $K_\circ |^\mathcal{B}$ is defined the same way as $K_\circ$ but with codomain $\mathcal{B}$.
\begin{theorem} \label{thm5.23}
    Suppose $\mathcal{D}$ is small, and $\mathcal{C}$ is locally small and complete. Let $(K,\kappa,\iota) \colon$ $\mathcal{D} \to \mathcal{C}$ be a lax monoidal functor, and $(\mathbb{T},\varepsilon)$ be the codensity monad of $K$. Then $(\mathbb{T},\varepsilon)$ is exactly pointwise monoidal iff there's a full monoidal subcategory $\mathcal{E}$ of $[\mathcal{D},\mathbf{Set}]^{\mathrm{op}}$, and an oplax monoidal coreflective subcategory $\mathcal{B}$ of $\mathcal{E}$ containing $\mathcal{A}$
    \[\begin{tikzcd}
	{\mathcal{C}} & {\mathcal{B}} & {\mathcal{E}}
	\arrow[""{name=0, anchor=center, inner sep=0}, "{K_\circ|^\mathcal{B}}", shift left=2, from=1-1, to=1-2]
	\arrow[""{name=1, anchor=center, inner sep=0}, "{\mathrm{Ran}_{Y_\mathcal{D}}K|_\mathcal{B}}", shift left=2, from=1-2, to=1-1]
	\arrow[""{name=2, anchor=center, inner sep=0}, "i", shift left=2, from=1-2, to=1-3]
	\arrow[""{name=3, anchor=center, inner sep=0}, "L", shift left=2, from=1-3, to=1-2]
	\arrow["\dashv"{anchor=center, rotate=-90}, draw=none, from=0, to=1]
	\arrow["\dashv"{anchor=center, rotate=-90}, draw=none, from=2, to=3]
\end{tikzcd}\]
such that
\begin{enumerate}
    \item The adjunction $K_\circ|^\mathcal{B} \dashv \mathrm{Ran}_{Y_\mathcal{D}}K |_\mathcal{B}$ lifts to $\mathbf{MonCat}$
    \item $iK_\circ|^\mathcal{B}=K_\circ|^\mathcal{E}$ as oplax monoidal functors and $\mathrm{Ran}_{Y_\mathcal{D}} K|_\mathcal{B}L=\mathrm{Ran}_{Y_\mathcal{D}} K|_\mathcal{E}$ as lax monoidal functors
\end{enumerate}
In particular, the Kleisli category of $\mathbb{T}$ is monoidally equivalent to an oplax monoidal coreflective subcategory of a monoidal subcategory of $[\mathcal{D},\mathbf{Set}]^{\mathrm{op}}$.
\end{theorem}
 \proof
 First, suppose $(\mathbb{T},\varepsilon)$ is exactly pointwise monoidal, and let $\Bar{\mathcal{A}}$ be the smallest full monoidal subcategory of $[\mathcal{D},\textbf{Set}]^{\text{op}}$ containing $\mathcal{A}$. Explicitly, objects of $\Bar{\mathcal{A}}$ are of the form $\bigotimes^\Lambda \mathcal{C}(A_i,K{-})$  where $\Lambda$ encodes a binary bracketing of $n\geq 0$ letters and special symbol $I$, and $\bigotimes^\Lambda X_i$ denotes substitution of the objects $X_1,\dots,X_n$ into $\Lambda$, for the binary operation $\bigotimes$. 
 
 Now, if $F \colon \mathcal{D} \to \textbf{Set}$ then a natural transformation $F \to \mathcal{C}(A,K{-})$ corresponds to a cone with apex $A$ over $KU \colon (* \downarrow F) \to \mathcal{C}$, and thus there is an isomorphism $[\mathcal{D},\textbf{Set}](F,\mathcal{C}(A,K{-})) \cong \mathcal{C}(A,\lim_{(* \downarrow F)}KU)$. In particular, by \hyperref[lem5.21]{Lemma 5.21}, for $k\geq 0$, $[\mathcal{D},\textbf{Set}](\bigotimes_{i=1}^k\mathcal{C}(A_i,K{-}),\mathcal{C}(B,K{-})) \cong \mathcal{C}(B, T(\bigotimes_{i=1}^k A_i))$. Hence, by coherence of monoidal categories ~\cite{mac_lane_categories_1998}, we see that $\mathcal{A}^\text{op}$ is a reflective subcategory of $\Bar{\mathcal{A}}^\text{op}$ where the unit of the reflection has component $\xi^{(k)}$ at $\bigotimes_{i=1}^k \mathcal{C}(A_i,K{-})$. This reflection satisfies the conditions of \hyperref[prop5.22]{Proposition 5.22}, so $\mathcal{A}$ is an oplax monoidal coreflective subcategory of $\Bar{\mathcal{A}}$ as desired. One can then verify that $K_\circ |^\mathcal{\bar{A}}$ is strong monoidal, and so by \hyperref[prop5.20]{Proposition 5.20} this satisfies the desired conditions. 
 
 Conversely, suppose we have $\mathcal{B}$, $\mathcal{E}$, and $L$ satisfying the conditions of the theorem. By \hyperref[prop5.20]{Proposition 5.20}, $K_\circ|^\mathcal{B}$ is strong monoidal. Hence, the maps $\bigotimes_{i=1}^k \mathcal{C}(A_i, K{-}) \to \mathcal{C}(\bigotimes_{i=1}^k A_i, K{-})$ induced by the oplax monoidal structure of $i$ are given by $\xi^{(k)}$. By \hyperref[prop5.22]{Proposition 5.22}, $L\xi^{(k)}$ are isomorphisms for $k\geq 0$. Hence, $\mathrm{Ran}_{Y_\mathcal{D}}K(\xi^{(k)})$ are also isomorphisms, and so by \hyperref[lem5.21]{Lemma 5.21}, $(\mathbb{T},\varepsilon)$ is exactly pointwise monoidal.
 \endproof
 \begin{example}\textbf{Giry monad on \textbf{BorelMeas}.}

We saw in \hyperref[example5.16]{Example 5.16} that $\mathcal{G}_B$ is exactly pointwise monoidal. \textbf{BorelMeas} is not complete, so we cannot immediately apply \hyperref[thm5.23]{Theorem 5.23}, but the proof above still shows that the Kleisli category of $\mathcal{G}_B$ is an oplax monoidal coreflective subcategory of $\Bar{\mathcal{A}}$, where $\Bar{\mathcal{A}}$ is the submonoidal category of $[\textbf{cStoch},\textbf{Set}]^\text{op}$ defined in the proof of \hyperref[thm5.23]{Theorem 5.23}.
 \end{example}
 An exactly pointwise monoidal codensity monad $(\mathbb{T},\varepsilon)$ of a lax monoidal functor $K$ is said to be exactly pointwise symmetric monoidal if $K$ is symmetric monoidal. Then \hyperref[prop5.20]{Proposition 5.20}, \hyperref[lem5.21]{Lemma 5.21}, \hyperref[prop5.22]{Proposition 5.22} and \hyperref[thm5.23]{Theorem 5.23}  remain true if you substitute all instances of \textit{monoidal} for \textit{symmetric monoidal}~\cite[Theorem 5.48]{malkiewich_coherence_2022}.
\begin{remark}
There are several definitions for when a Kan extension in a 2{-}category is pointwise. \cite{street_fibrations_1974} provides one for a 2{-}category with \emph{comma objects}. These are lax squares of the form
\[\begin{tikzcd}
	{(H\downarrow F)} & {\mathcal{E}} \\
	{\mathcal{C}} & {\mathcal{D}}
	\arrow["{U_H}", from=1-1, to=1-2]
	\arrow["{U_F}"', from=1-1, to=2-1]
	\arrow["\alpha"', shorten <=10pt, shorten >=10pt, Rightarrow, from=1-2, to=2-1]
	\arrow["H", from=1-2, to=2-2]
	\arrow["F"', from=2-1, to=2-2]
\end{tikzcd}\]
that have a universal property that can be found in \cite{street_fibrations_1974}. We refer to this square as a comma square. In this setting, a right Kan extension along $F$ is defined to be pointwise if, for any $H$, it is preserved by pasting with this comma square. For example, if $H,F$ are functors, then the comma category $(H\downarrow F)$ is a comma object in $\textbf{Cat}$, and this definition recovers the notion of a pointwise right Kan extension in \textbf{Cat}. When $(H,\chi,\nu)$ is strong monoidal, and $(F,\kappa,\iota)$ is lax monoidal, we can give a monoidal category structure to $(H\downarrow F)$ defined on objects by $(A,X,h)\otimes (B,Y,k)=(A\otimes B, X \otimes Y, \kappa_{X,Y}(h\otimes k)\chi_{A,B}^{{-}1})$ with unit $(I,I,\iota \nu^{{-}1})$. Then $U_H, U_F$ are strict monoidal, $\alpha$ is monoidal, and this provides a comma object in $\textbf{MonCat}$. Now, if $\text{Ran}_F G$ is an exactly pointwise monoidal Kan extension (is preserved by pasting the lax morphism square for $F$), one can verify that $\text{Ran}_F G$ is preserved by pasting comma squares in $\textbf{MonCat}$ for any strong monoidal $H$, and that the resulting Kan extension is exactly pointwise monoidal. However, it is not clear that these are the only Kan extensions preserved by pasting these comma squares, and whether these are the only comma objects that exist in $\textbf{MonCat}$.
\end{remark}
\bibliography{references}

@inproceedings{jones_probabilistic_1989,
	title = {A probabilistic powerdomain of evaluations},
	isbn = {978-0-8186-1954-0},
	url = {https://ieeexplore.ieee.org/document/39173},
	doi = {10.1109/LICS.1989.39173},
	language = {eng},
	urldate = {2025-07-31},
	booktitle = {Proc. {LICS}},
	author = {Jones, C. and Plotkin, G. D.},
	year = {1989},
	pages = {186--195},
}

@misc{van_breugel_metric_nodate,
	title = {The {Metric} {Monad} for {Probabilistic} {Nondeterminism}},
	url = {https://citeseerx.ist.psu.edu/document?repid=rep1&type=pdf&doi=2c172448c5ea3504bcae6b49ed47eb113eca51ef},
	urldate = {2023-10-13},
	author = {van Breugel, Franck},
	note = {Draft available at http://www.cse.yorku.ca/\%7Efranck/research/drafts/monad.pdf, 2005},
}

@misc{fritz_criterion_nodate,
	title = {A {Criterion} for {Kan} {Extensions} of {Lax} {Monoidal} {Functors}},
	doi = {10.48550/arxiv.1809.10481},
	language = {eng},
	author = {Fritz, Tobias and Perrone, Paolo},
	note = {arxiv:1809.10481, 2018},
}

@misc{fritz_convex_nodate,
	title = {Convex {Spaces} {I}: {Definition} and {Examples}},
	shorttitle = {Convex {Spaces} {I}},
	doi = {10.48550/arxiv.0903.5522},
	language = {eng},
	author = {Fritz, Tobias},
	note = {arXiv:0903.5522, 2009},
}

@misc{basso_hitchhikers_2015,
	title = {A {Hitchhiker}’s guide to {Wasserstein} distances},
	url = {https://www.semanticscholar.org/paper/A-Hitchhiker-%E2%80%99-s-guide-to-Wasserstein-distances-Basso/619460ff6f03ecb4c19b91b9c195c5bfac3e8f06},
	urldate = {2024-04-01},
	author = {Basso, Giuliano},
	year = {2015},
}

@article{kelly_doctrinal_1974,
	title = {Doctrinal {Adjunction}},
	volume = {420},
	copyright = {http://www.springer.com/tdm},
	journal = {Lecture Notes in Mathematics},
	author = {Kelly, G. M.},
	year = {1974},
	pages = {257--280},
}

@phdthesis{koudenburg_algebraic_2013,
	title = {Algebraic weighted colimits},
	url = {http://arxiv.org/abs/1304.4079},
	urldate = {2025-08-21},
	school = {University of Sheffield},
	author = {Koudenburg, Seerp Roald},
	year = {2013},
	note = {arXiv:1304.4079},
}

@book{manes_algebraic_1976,
	address = {New York},
	series = {Graduate {Texts} in {Mathematics}},
	title = {Algebraic theories},
	volume = {26},
	isbn = {978-0-387-90140-4},
	language = {eng},
	publisher = {Springer-Verlag},
	author = {Manes, Ernest G.},
	year = {1976},
}

@article{kennison_equational_1971,
	title = {Equational completion, model induced triples and pro-objects},
	volume = {1},
	issn = {0022-4049},
	doi = {10.1016/0022-4049(71)90001-6},
	language = {eng},
	number = {4},
	journal = {Journal of Pure and Applied Algebra},
	author = {Kennison, J. F. and Gildenhuys, Dion},
	year = {1971},
	pages = {317--346},
}

@article{moggi_notions_1991,
	title = {Notions of computation and monads},
	volume = {93},
	issn = {0890-5401},
	doi = {10.1016/0890-5401(91)90052-4},
	language = {eng},
	number = {1},
	journal = {Information and computation},
	author = {Moggi, Eugenio},
	year = {1991},
	pages = {55--92},
}

@misc{lawvere_category_1962,
	title = {The category of probabilistic mappings},
	url = {https://lawverearchives.com/wp-content/uploads/2025/07/1962.probmap.pdf},
	urldate = {2025-08-14},
	author = {Lawvere, F. William},
	year = {1962},
	note = {Online pdf},
}

@inproceedings{heunen_convenient_2017,
	title = {A convenient category for higher-order probability theory},
	isbn = {978-1-5090-3019-4},
	booktitle = {32nd {Annual} {ACM}/{IEEE} {Symposium} on {Logic} in {Computer} {Science}},
	author = {Heunen, Chris and Kammar, Ohad and Staton, Sam and Yang, Hongseok},
	year = {2017},
}

@article{fritz_finettis_2021,
	title = {De {Finetti}’s theorem in categorical probability},
	volume = {2},
	issn = {2689-6931},
	language = {eng},
	number = {4},
	journal = {Journal of Stochastic Analysis},
	author = {Fritz, Tobias and Gonda, Tomáš and Perrone, Paolo},
	year = {2021},
}

@article{epstein_functors_1966,
	title = {Functors between {Tensored} {Categories}.},
	volume = {1},
	issn = {0020-9910; 1432-1297/e},
	url = {https://eudml.org/doc/141834},
	language = {und},
	urldate = {2025-05-06},
	journal = {Inventiones mathematicae},
	author = {Epstein, D. B. A.},
	year = {1966},
	pages = {221--228},
}

@article{jamneshan_uncountable_2023,
	title = {An uncountable {Moore}–{Schmidt} theorem},
	volume = {43},
	issn = {0143-3857},
	doi = {10.1017/etds.2022.36},
	language = {eng},
	number = {7},
	journal = {Ergodic theory and dynamical systems},
	author = {Jamneshan, Asgar and Tao, Terrence},
	year = {2023},
	pages = {2376--2403},
}

@incollection{mulry_lifting_1993,
	address = {Berlin},
	series = {Lecture {Notes} in {Computer} {Science}},
	title = {Lifting theorems for {Kleisli} categories},
	volume = {802},
	isbn = {978-3-540-58027-0},
	language = {eng},
	booktitle = {Mathematical {Foundations} of {Programming} {Semantics}},
	publisher = {Springer},
	author = {Mulry, Philip S.},
	editor = {Melton, Austin and Brookes, Stephen and Main, Michael and Mislove, Michael and Schmidt, David},
	year = {1993},
	pages = {304--319},
}

@article{malkiewich_coherence_2022,
	title = {Coherence for bicategories, lax functors, and shadows},
	volume = {38},
	issn = {1201-561X},
	language = {eng},
	number = {12},
	journal = {Theory and Applications of Categories},
	author = {Malkiewich, Cary and Ponto, Kate},
	year = {2022},
	pages = {328--373},
}

@article{tadeusz_swirszcz_monadic_1974,
	title = {Monadic {Functors} and {Convexity}},
	volume = {22},
	journal = {Bulletin de l’Académie Polonaise des Sciences: Série des Sciences Mathématiques, Astronomiques et Physiques},
	author = {{Tadeusz Świrszcz}},
	year = {1974},
}

@article{cho_disintegration_2019,
	title = {Disintegration and {Bayesian} {Inversion} via {String} {Diagrams}},
	volume = {29},
	issn = {0960-1295},
	doi = {10.1017/S0960129518000488},
	language = {eng},
	number = {7},
	journal = {Mathematical Structures in Computer Science},
	author = {Cho, Kenta and Jacobs, Bart},
	year = {2019},
	pages = {938--971},
}

@book{huber_robust_1981,
	address = {New York},
	series = {Wiley {Series} in {Probability} and {Statistics}},
	title = {Robust {Statistics}},
	isbn = {978-0-471-41805-4},
	language = {eng},
	publisher = {Wiley},
	author = {Huber, Peter J.},
	year = {1981},
}

@book{billingsley_convergence_1999,
	address = {New York},
	series = {Wiley {Series} in {Probability} and {Statistics}},
	title = {Convergence of {Probability} {Measures}},
	isbn = {978-1-118-62596-5},
	language = {eng},
	publisher = {Wiley},
	author = {Billingsley, Patrick},
	year = {1999},
}

@book{dudley_real_2002,
	address = {Cambridge},
	series = {Cambridge {Studies} in {Advanced} {Mathematics}},
	title = {Real {Analysis} and {Probability}},
	isbn = {978-1-351-09309-5},
	language = {eng},
	publisher = {Cambridge University Press},
	author = {Dudley, Richard M.},
	year = {2002},
	doi = {10.1201/9781351076197},
}

@book{parthasarathy_probability_1967,
	address = {New York, London},
	series = {Probability and {Mathematical} {Statistics}: {A} {Series} of {Monographs} and {Textbooks}},
	title = {Probability {Measures} on {Metric} {Spaces}},
	isbn = {978-0-8218-3889-1},
	language = {eng},
	publisher = {Academic Press},
	author = {Parthasarathy, Kalyanapuram R.},
	year = {1967},
}

@book{mac_lane_categories_1998,
	address = {New York},
	series = {Graduate {Texts} in {Mathematics}},
	title = {Categories for the {Working} {Mathematician}},
	isbn = {978-1-4419-3123-8},
	language = {eng},
	publisher = {Springer},
	author = {Mac Lane, Saunders},
	year = {1998},
}

@book{barr_toposes_1985,
	address = {New York},
	series = {Grundlehren der mathematischen {Wissenschaften}},
	title = {Toposes, {Triples} and {Theories}},
	isbn = {978-0-387-96115-6},
	language = {eng},
	publisher = {Springer},
	author = {Barr, Michael},
	collaborator = {Wells, Charles},
	year = {1985},
}

@article{fritz_bimonoidal_2018,
	title = {Bimonoidal {Structure} of {Probability} {Monads}},
	volume = {341},
	issn = {1571-0661},
	doi = {10.1016/j.entcs.2018.11.007},
	language = {eng},
	journal = {Electronic Notes in Theoretical Computer Science},
	author = {Fritz, Tobias and Perrone, Paolo},
	year = {2018},
	pages = {121--149},
}

@incollection{street_fibrations_1974,
	address = {Berlin},
	series = {Lecture {Notes} in {Mathematics}},
	title = {Fibrations and {Yoneda}'s lemma in a 2-category},
	volume = {420},
	isbn = {978-3-540-06966-9},
	language = {eng},
	booktitle = {Category {Seminar}: {Proceedings} of the {Sydney} {Category} {Theory} {Seminar}, 1972/1973},
	publisher = {Springer},
	author = {Street, Ross},
	editor = {{Kelly, G.}},
	year = {1974},
	pages = {104--133},
}

@incollection{day_closed_1970,
	address = {Berlin},
	series = {Lecture {Notes} in {Mathematics}},
	title = {On closed categories of functors},
	volume = {137},
	isbn = {978-3-540-04926-5},
	language = {eng},
	booktitle = {Reports of the {Midwest} {Category} {Seminar} {IV}},
	publisher = {Springer},
	author = {Day, Brian},
	editor = {Mac Lane, S.},
	year = {1970},
	pages = {1--38},
}

@article{van_belle_probability_2022,
	title = {Probability monads as codensity monads},
	volume = {38},
	issn = {1201-561X},
	language = {eng},
	number = {21},
	journal = {Theory and Applications of Categories},
	author = {Van Belle, Ruben},
	year = {2022},
	pages = {811--842},
}

@book{bogachev_measure_2007,
	address = {Berlin},
	title = {Measure {Theory}},
	isbn = {978-3-540-34513-8},
	language = {eng},
	publisher = {Springer},
	author = {Bogachev, Vladimir I.},
	year = {2007},
	doi = {10.1007/978-3-540-34514-5},
}

@book{loregian_coend_2021,
	address = {Cambridge},
	series = {London {Mathematical} {Society} {Lecture} {Note} {Series}},
	title = {({Co})end calculus},
	isbn = {978-1-108-78860-1},
	language = {eng},
	publisher = {Cambridge University Press},
	author = {Loregian, Fosco},
	year = {2021},
}

@article{day_note_1973,
	title = {Note on monoidal localisation},
	volume = {8},
	issn = {0004-9727},
	doi = {10.1017/S0004972700045433},
	language = {eng},
	number = {1},
	journal = {Bulletin of the Australian Mathematical Society},
	author = {Day, Brian},
	year = {1973},
	pages = {1--16},
}

@article{im_universal_1986,
	title = {A universal property of the convolution monoidal structure},
	volume = {43},
	issn = {0022-4049},
	doi = {10.1016/0022-4049(86)90005-8},
	language = {eng},
	number = {1},
	journal = {Journal of Pure and Applied Algebra},
	author = {Im, Geun Bin and Kelly, Gregory Maxwell},
	year = {1986},
	pages = {75--88},
}

@article{koudenburg_algebraic_2015,
	title = {Algebraic {Kan} extensions in double categories},
	volume = {30},
	issn = {1201-561X},
	language = {eng},
	number = {5},
	journal = {Theory and Applications of Categories},
	author = {Koudenburg, Seerp Roald},
	year = {2015},
	pages = {86--146},
}

@article{weber_algebraic_2016,
	title = {Algebraic {Kan} extensions along morphisms of internal algebra classifiers},
	volume = {9},
	issn = {1875-158X},
	doi = {10.1515/tmj-2016-0006},
	language = {eng},
	number = {1},
	journal = {Tbilisi Mathematical Journal},
	author = {Weber, Mark},
	year = {2016},
	pages = {65--142},
}

@book{semadeni_monads_1973,
	address = {Kingston, Ontario},
	series = {Queen's {Papers} in {Pure} and {Applied} {Mathematics}},
	title = {Monads and their {Eilenberg}-{Moore} {Algebras} in {Functional} {Analysis}},
	language = {eng},
	publisher = {Queen's University},
	author = {Semadeni, Zbigniew},
	year = {1973},
}

@book{halmos_measure_1974,
	address = {New York},
	series = {Graduate {Texts} in {Mathematics}},
	title = {Measure theory},
	isbn = {978-0-387-90088-9},
	language = {eng},
	publisher = {Springer},
	author = {Halmos, Paul R.},
	year = {1974},
}

@article{ulam_zur_1930,
	title = {Zur {Masstheorie} in der allgemeinen {Mengenlehre}},
	volume = {16},
	issn = {0016-2736},
	doi = {10.4064/fm-16-1-140-150},
	language = {eng},
	journal = {Fundamenta Mathematicae},
	author = {Ulam, Stanisław},
	year = {1930},
	pages = {140--150},
}

@book{rudin_functional_1991,
	address = {New York},
	series = {International {Series} in {Pure} and {Applied} {Mathematics}},
	title = {Functional {Analysis}},
	isbn = {978-0-07-054236-5},
	language = {eng},
	publisher = {McGraw-Hill},
	author = {Rudin, Walter},
	year = {1991},
}

@article{marshall_h_stone_postulates_1949,
	title = {Postulates for the barycentric calculus},
	volume = {29},
	issn = {0373-3114},
	doi = {10.1007/bf02413910},
	language = {eng},
	journal = {Annali di Matematica Pura ed Applicata},
	author = {{Marshall H. Stone}},
	year = {1949},
	pages = {25--30},
}

@article{leinster_codensity_2013,
	title = {Codensity and the ultrafilter monad},
	volume = {28},
	issn = {1201-561X},
	language = {eng},
	number = {1},
	journal = {Theory and Applications of Categories},
	author = {Leinster, Tom},
	year = {2013},
	pages = {332--370},
}

@article{kock_strong_1972,
	title = {Strong functors and monoidal monads},
	volume = {23},
	issn = {0003-889X},
	doi = {10.1007/BF01304852},
	language = {eng},
	number = {1},
	journal = {Archiv der Mathematik},
	author = {Kock, Anders},
	year = {1972},
	pages = {113--120},
}

@article{kock_closed_1971,
	title = {Closed categories generated by commutative monads},
	volume = {12},
	issn = {0004-9735},
	doi = {10.1017/S1446788700010272},
	language = {eng},
	number = {4},
	journal = {Journal of the Australian Mathematical Society},
	author = {Kock, Anders},
	year = {1971},
	pages = {405--424},
}

@article{kock_monads_1970,
	title = {Monads on symmetric monoidal closed categories},
	volume = {21},
	issn = {0003-889X},
	doi = {10.1007/BF01220868},
	language = {eng},
	number = {1},
	journal = {Archiv der Mathematik},
	author = {Kock, Anders},
	year = {1970},
	pages = {1--10},
}

@article{keimel_monad_2008,
	title = {The monad of probability measures over compact ordered spaces and its {Eilenberg}–{Moore} algebras},
	volume = {156},
	issn = {0166-8641},
	doi = {10.1016/j.topol.2008.07.002},
	language = {eng},
	number = {2},
	journal = {Topology and its Applications},
	author = {Keimel, Klaus},
	year = {2008},
	pages = {227--239},
}

@article{karni_extension_1990,
	title = {On the extension of bimeasures},
	volume = {55},
	issn = {1565-8538},
	url = {https://doi.org/10.1007/BF02789194},
	doi = {10.1007/BF02789194},
	language = {en},
	number = {1},
	urldate = {2024-04-07},
	journal = {Journal d’Analyse Mathématique},
	author = {Karni, Shaul and Merzbach, Ely},
	year = {1990},
	pages = {1--16},
}

@article{jacobs_probability_2018,
	title = {From probability monads to commutative effectuses},
	volume = {94},
	issn = {2352-2208},
	doi = {10.1016/j.jlamp.2016.11.006},
	language = {eng},
	journal = {Journal of Logical and Algebraic Methods in Programming},
	author = {Jacobs, Bart},
	year = {2018},
	pages = {200--237},
}

@incollection{giry_categorical_1982,
	address = {Berlin},
	series = {Lecture {Notes} in {Mathematics}},
	title = {A categorical approach to probability theory},
	volume = {915},
	isbn = {978-3-540-11211-2},
	language = {eng},
	booktitle = {Categorical {Aspects} of {Topology} and {Analysis}},
	publisher = {Springer},
	author = {Giry, Michèle},
	editor = {Banaschewski, B.},
	year = {1982},
	pages = {68--85},
}

@article{fritz_probability_2019,
	title = {A {Probability} {Monad} as the {Colimit} of {Spaces} of {Finite} {Samples}},
	volume = {34},
	issn = {1201-561X},
	language = {eng},
	number = {7},
	journal = {Theory and Applications of Categories},
	author = {Fritz, Tobias and Perrone, Paolo},
	year = {2019},
	pages = {170--220},
}

@article{fritz_synthetic_2020,
	title = {A synthetic approach to {Markov} kernels, conditional independence and theorems on sufficient statistics},
	volume = {370},
	issn = {0001-8708},
	doi = {10.1016/j.aim.2020.107239},
	language = {eng},
	journal = {Advances in Mathematics},
	author = {Fritz, Tobias},
	year = {2020},
	pages = {107239},
}

@book{fremlin_measure_2000,
	address = {Colchester},
	title = {Measure {Theory}},
	isbn = {978-0-9538129-0-5},
	language = {eng},
	publisher = {Torres Fremlin},
	author = {Fremlin, David H.},
	year = {2000},
}

@article{fedorchuk_probability_1991,
	title = {Probability measures in topology},
	volume = {46},
	issn = {0036-0279},
	doi = {10.1070/RM1991v046n01ABEH002722},
	language = {eng},
	number = {1},
	journal = {Russian Mathematical Surveys},
	author = {Fedorchuk, Vitalii V.},
	year = {1991},
	pages = {45--93},
}

@article{dobrakov_multilinear_1999,
	title = {Multilinear integration of bounded scalar valued functions},
	volume = {49},
	issn = {0232-0525},
	url = {https://eudml.org/doc/34496},
	language = {eng},
	number = {3},
	urldate = {2024-04-07},
	journal = {Mathematica Slovaca},
	author = {Dobrakov, Ivan},
	year = {1999},
	pages = {295--304},
}

@article{bowers_representation_2015,
	title = {Representation of {Extendible} {Bilinear} {Forms}},
	volume = {65},
	issn = {0139-9918},
	doi = {10.1515/ms-2015-0077},
	language = {eng},
	number = {5},
	journal = {Mathematica Slovaca},
	author = {Bowers, Adam},
	year = {2015},
	pages = {1123--1136},
}

@article{bombal_integral_2001,
	title = {Integral {Operators} on the {Product} of {C}({K}) {Spaces}},
	volume = {264},
	issn = {0022-247X},
	doi = {10.1006/jmaa.2001.7648},
	language = {eng},
	number = {1},
	journal = {Journal of Mathematical Analysis and Applications},
	author = {Bombal, Fernando and Villanueva, Ignacio},
	year = {2001},
	pages = {107--121},
}

@article{avery_codensity_2016,
	title = {Codensity and the {Giry} monad},
	volume = {220},
	issn = {0022-4049},
	doi = {10.1016/j.jpaa.2015.08.017},
	language = {eng},
	number = {3},
	journal = {Journal of Pure and Applied Algebra},
	author = {Avery, Tom},
	year = {2016},
	pages = {1229--1251},
}

@article{marczewski_remarks_1953,
	title = {Remarks on the compactness and non direct products of measures},
	volume = {40},
	issn = {0016-2736},
	url = {https://eudml.org/doc/213318},
	language = {eng},
	number = {1},
	urldate = {2024-04-07},
	journal = {Fundamenta Mathematicae},
	author = {Marczewski, Edward and Ryll-Nardzewski, Czesław},
	year = {1953},
	pages = {165--170},
}

@book{halmos_lectures_1974,
	address = {New York},
	title = {Lectures on {Boolean} algebras},
	isbn = {978-0-387-90094-0},
	language = {eng},
	publisher = {Springer-Verlag},
	author = {Halmos, Paul R.},
	year = {1974},
}
\appendix
\section{Probability measures on compact Hausdorff spaces}
\label{appenA}
 In this appendix, we review some basic theory of Baire measures and prove some useful properties of the functor $H \colon \textbf{KHaus} \to \textbf{Meas}$, which assigns compact Hausdorff spaces their Baire measurable structure.
\begin{definition}[Baire $\sigma${-}algebra~\cite{halmos_measure_1974}] If $X$ is a topological space, its Baire $\sigma${-}algebra $\mathcal{H}X$ is the coarsest $\sigma$ algebra such that each continuous $h\colon X \to [0,1]$ is also measurable. 
\end{definition}
  For a compact Hausdorff space $X$, $\mathcal{H}X$ is equivalently given by the $\sigma${-}algebra generated by the zero sets in $X$, which are sets of the form $f^{{-}1}(\{0\})$ for continuous $f \colon X \to [0,1]$. For any metric space $X$, we have $\mathcal{H}X=\mathcal{B}X$~\cite[Theorem 7.1.1]{dudley_real_2002}. However for general topological spaces $X$, we only have $\mathcal{H}X \subseteq \mathcal{B}X$, and this inclusion can be strict. For $X=\{0,1\}^\kappa$, when $\kappa >\omega$ the points of $X$ are closed and are Borel measurable, but by \hyperref[thmA.3]{Theorem A.3} they are not Baire measurable. In fact, for a Stone space $X$, $\mathcal{H}X$ is the $\sigma${-}algebra generated by its clopen sets~\cite[Section 51]{halmos_measure_1974}.
  
  For the remainder of this section, $X$ will be a compact Hausdorff space. Let $B \colon \allowbreak \textbf{KHaus} \to \textbf{Meas}$ be the functor assigning a space its Borel $\sigma${-}algebra, and let $H$ assign the Baire $\sigma${-}algebra as above (both are identity on morphisms). There is a natural map $B \to H$ given by the inclusion $\mathcal{H}X \subseteq \mathcal{B}X$. This induces a measurable natural map $\mathcal{G}B \to \mathcal{G}H$ which restricts Borel probability measures to Baire probability measures. For the subspace of Radon measures in $\mathcal{G}BX$, this map restricts to a bijection.
 
\begin{proposition}[Theorem 7.3.1~\cite{dudley_real_2002}] \label{propA.2} If $X$ is a compact Hausdorff space and $p$ is a probability measure on $\mathcal{H}X$ then there is a unique Radon measure on $\mathcal{B}X$ restricting to $p$.
\end{proposition}

 Now, let $\mathcal{A}$ be the coarsest $\sigma${-}algebra on the underlying set of $\mathcal{G}HX$ such that \allowbreak $\text{ev}_h \colon \mathcal{G}HX \to \mathbb{R}$ is measurable for each continuous $h \colon X \to \mathbb{R}$. If $\mathcal{B}$ is the standard $\sigma${-}algebra on $\mathcal{G}HX$, then $\mathcal{A} \subseteq \mathcal{B}$. Let $V$ be the vector space of measurable functions $m \colon HX \to \mathbb{R}$ such that $\text{ev}_m \colon (\mathcal{G}HX,\mathcal{A}) \to \mathbb{R}$ is measurable. Then $C(X) \leq V$, and by the dominated convergence theorem, $V$ is closed under bounded pointwise limits, so $V$ contains all bounded measurable functions $HX \to \mathbb{R}$, and hence $\mathcal{A}=\mathcal{B}$. The function $\mathrm{res}_X \colon  H\mathcal{R}X \to \mathcal{G}HX$ that restricts Radon measures on $\mathcal{B}X$ to $\mathcal{H}X$ is therefore measurable. By \hyperref[propA.2]{Proposition A.2}, $\mathrm{res}_X$ is a bijection. To show it is an isomorphism, we use the following characterisation of the Baire $\sigma${-}algebra, which generalises results in \cite{halmos_measure_1974,jamneshan_uncountable_2023}.
\begin{theorem} \label{thmA.3}
    $H \colon \mathbf{KHaus} \to \mathbf{Meas}$ is the unique functor, up to isomorphism, that preserves limits and such that $H[0,1]\cong[0,1]$ (with its Borel $\sigma${-}algebra).
\end{theorem}
\proof
First, let $X$ be a compact Hausdorff space and $Y \subseteq X$ be a closed subspace. Let $\mathcal{H}Y$ be its Baire $\sigma${-}algebra, and $\mathcal{H}'Y$ be the subspace $\sigma${-}algebra of $Y \subseteq HX$. $\mathcal{H}'Y$ is the coarsest $\sigma${-}algebra such that $h|_Y \colon Y \to [0,1]$ is measurable for each continuous $h \colon X \to [0,1]$. Since each such $h|_Y$ is continuous, $\mathcal{H}'Y \subseteq \mathcal{H}Y$. Also, by the Tietze extension theorem, for any continuous $k \colon Y \to [0,1]$ there is an $h \colon X \to [0,1]$ so that $h|_Y=k$. Hence, $\mathcal{H}'Y=\mathcal{H}Y$ and $H$ preserves embeddings, so in particular preserves equalisers.

Now, let $X_i$ be compact Hausdorff spaces for $i \in I$ and let $X=\prod_{i \in I} X_i$. Let $\mathcal{H}X$ be its Baire $\sigma${-}algebra, and $\mathcal{H}'X$ be the $\sigma${-}algebra of the product $\prod_{i \in I} HX_i$ in \textbf{Meas}. By definition, $\mathcal{H}'X \subseteq \mathcal{H}X$, and for the reverse inclusion, we need to show that any continuous $m \colon X \to [0,1]$ is measurable for $\mathcal{H}'X$. Now, consider the subalgebra of $C(X)$ generated by finite products of continuous maps in one variable. Explicitly these are $m$ of the form $m(\underline{x})=\prod_{j=1}^k h_{i_j}(x_{i_j})$ where $h_{i_j} \colon X_{i_j} \to \mathbb{R}$ is continuous. Clearly each such $m$ is measurable for $\mathcal{H}'X$, and by the Stone{-}Weierstrass theorem, this is a dense subalgebra in $C(X)$, so it follows that $\mathcal{H}'X = \mathcal{H}X$.

Conversely, suppose $H'$ is any functor that preserves limits, and such that $H'[0,1]\cong[0,1]$.  For any compact Hausdorff space $X$, there is an embedding $X \to [0,1]^{C(X,[0,1])}$. Since all subspace inclusions are regular monic in $\textbf{KHaus}$, it follows that $H'\cong H$. 
\endproof
\begin{corollary}
\label{corA.4}
     The map $\mathrm{res} \colon H\mathcal{R} \to \mathcal{G}H$ is an isomorphism. Also if we let  $I \colon \allowbreak \mathbf{FinStoch} \to \mathbf{cStoch}$ be the inclusion, then we have $\mathcal{G}H=\textnormal{Ran}_{K_\mathcal{R}}K_\mathcal{G}I$.
\end{corollary}
\proof
    Since $\mathcal{R}=\text{Ran}_{K_\mathcal{R}} K_\mathcal{R}$ is pointwise, $H$ preserves limits and $HK_\mathcal{R}= K_\mathcal{G}I$, it follows that  $H\mathcal{R}=\text{Ran}_{K_\mathcal{R}}K_\mathcal{G}I$. Thus, there is a unique map $
    \alpha \colon \mathcal{G}H \to H\mathcal{R}$ satisfying $ \int h\; \text{d}\alpha_X(p)=\int h \; \text{d}p$ for each continuous $h \colon X \to [0,1]$. But this determines $p$ completely, and $\alpha$ is an inverse to $\text{res}$. 
\endproof
 As an illustrative example of the strength of this result, consider the functor $\mathcal{G}_f$ of finitely additive probability measures~\cite{van_belle_probability_2022} on \textbf{Meas}. There is a map $\mu^{\mathcal{G}_f} \colon \mathcal{G}_fK_\mathcal{G}I \to \mathcal{G}I$, and so \hyperref[corA.4]{Corollary A.4} produces a unique factorisation $\mathcal{G}_f H \to \mathcal{G}H$. This shows that the space $\mathcal{G}HX$ is a retract of $\mathcal{G}_f HX$. Hence, for any finitely additive Baire probability measure $p^*$, there is a countably additive Baire probability measure $p$ such that $\int h \; \text{d}p^*=\int h \; \text{d}p$ for every continuous $h \colon X \to [0,1]$.
\begin{example}
Let $X$ be a countable compact Hausdorff space, then both the Borel and Baire $\sigma${-}algebras of $X$ are the entire power set $\mathcal{P}X$. Now, let $\mathcal{U}$ be a non{-}principal ultrafilter on $X$, which can be regarded as a finitely additive measure.  Since $X$ is a compact Hausdorff space, $\mathcal{U}$ converges to a unique element $x\in X$. If $f \colon X \to [0,1]$ is continuous, then for every $\varepsilon >0$ we have $U=f^{{-}1}(f(x){-}\varepsilon,1]\in \mathcal{U}$ and the simple function $(f(x){-}\varepsilon) \mathbbm{1}_U \leq f$. Thus, $\int f \; \text{d}\, \mathcal{U} \geq f(x)$, and the reverse inequality is similar. Hence, we have  $\int f \; \text{d}\, \mathcal{U}=f(x)=\int f \; \text{d}\eta_X(x)$.
\end{example}
\end{document}